\def\jktr{n} \def\draft{n}

\documentclass[12pt]{amsart}
\usepackage{stmaryrd,fullpage,amsthm,amsmath,graphicx,verbatim,hyperref,picins}
\usepackage[all]{xy}
\usepackage{subcaption}

\if\draft y
  \usepackage[color]{showkeys}
  \usepackage{everypage,draftwatermark}
  \SetWatermarkScale{4}
  \SetWatermarkLightness{0.9}
\fi

\def\arXiv#1{{\href{http://front.math.ucdavis.edu/#1}{arXiv:#1}}}

\catcode`\@=11
\long\def\@makecaption#1#2{%
     \vskip 10pt
     \setbox\@tempboxa\hbox{
       \small\sf{\bfcaptionfont #1. }\ignorespaces #2}%
     \ifdim \wd\@tempboxa >\captionwidth {%
         \rightskip=\@captionmargin\leftskip=\@captionmargin
         \unhbox\@tempboxa\par}%
       \else
         \hbox to\hsize{\hfil\box\@tempboxa\hfil}%
     \fi}
\font\bfcaptionfont=cmssbx10 scaled \magstephalf
\newdimen\@captionmargin\@captionmargin=2\parindent
\newdimen\captionwidth\captionwidth=\hsize
\catcode`\@=12

\theoremstyle{plain}
\newtheorem{theorem}{Theorem}
\newtheorem{conjecture}{Conjecture}

\def\imagetop#1{\vtop{\null\hbox{#1}}}

\begin{document}

\title[Meta-Monoids]{Meta-Monoids, Meta-Bicrossed Products, and the Alexander Polynomial}

\author{Dror Bar-Natan}
\address{
  Department of Mathematics\\
  University of Toronto\\
  Toronto Ontario M5S 2E4\\
  Canada
}
\email{drorbn@math.toronto.edu}
\urladdr{http://www.math.toronto.edu/~drorbn}

\author{Sam Selmani}
\address{
  McGill University\\
  Department of Physics\\
  Montreal, Quebec H3A 2T8\\
  Canada
}
\email{sam.selmani@physics.mcgill.ca}

\date{\today; first edition: February 19, 2013. To appear in the {\em Journal of Knot Theory and its Ramifications}}

\thanks{This work was partially supported by NSERC grant RGPIN 262178 and partially pursued at the Newton Institute in Cambridge, UK. The full \TeX\ sources are at \url{http://drorbn.net/AcademicPensieve/Projects/MetaMonoids/}. This is \arXiv{1302.5689}.}

\maketitle

\begin{abstract}
We introduce a new invariant of tangles along with an algebraic framework in which to understand it. We claim that the invariant contains the classical Alexander polynomial of knots and its multivariable extension to links. We argue that of the computationally efficient members of the family of Alexander invariants, it is the most meaningful.

These are lecture notes for talks given by the first author, written and completed by the second. The talks, with handouts and videos, are available at \url{http://www.math.toronto.edu/drorbn/Talks/Regina-1206/}. See also further comments at \url{http://www.math.toronto.edu/drorbn/Talks/Caen-1206/#June8}.
\end{abstract}

\tableofcontents

\section{Warm-up: the baby invariant, $Z^G$}\label{one}

Let $T$ be an oriented tangle diagram. Let $G$ be a monoid\footnote{A monoid is like a group, but without inverses: it is a set with an associative binary operation and a unit. Every group is also a monoid.}, and suppose we are given two pairs $R^\pm=(g_o^\pm,g_u^\pm)$ of elements of $G$. At each positive (resp. negative)\footnotemark \space crossing of $T$, assign $g_o^+$ (resp. $g_o^-$) to the upper strand and $g_u^+$ (resp. $g_u^-$) to the lower strand, as in Figure \ref{fig:MonoidZ}. Then, for every strand, multiply all elements assigned to it in the order that they appear and store the end result. If $T$ has $n$ strands, we get a collection of $n$ elements of $G$. Call this collection $Z^G(T)$.
\footnotetext{Signs are determined by the ``right-hand rule'': If the right-hand thumb points along the direction of the upper strand of a positive crossing, then the fingers curl in the direction of the lower strand.}

\if\jktr y
  \begin{figure}[h]\centering
    \includegraphics[scale=0.54]{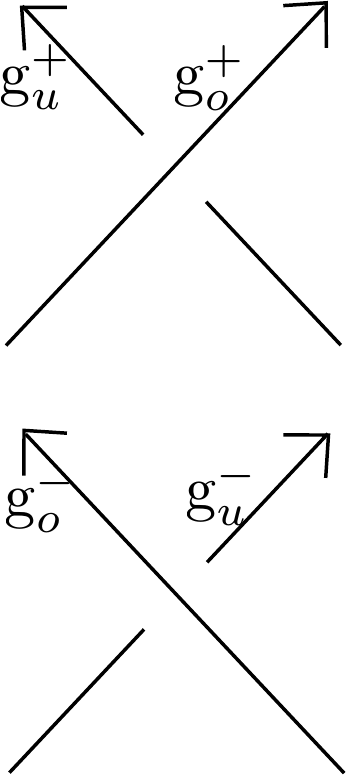}
    \qquad
    \includegraphics[scale=0.54]{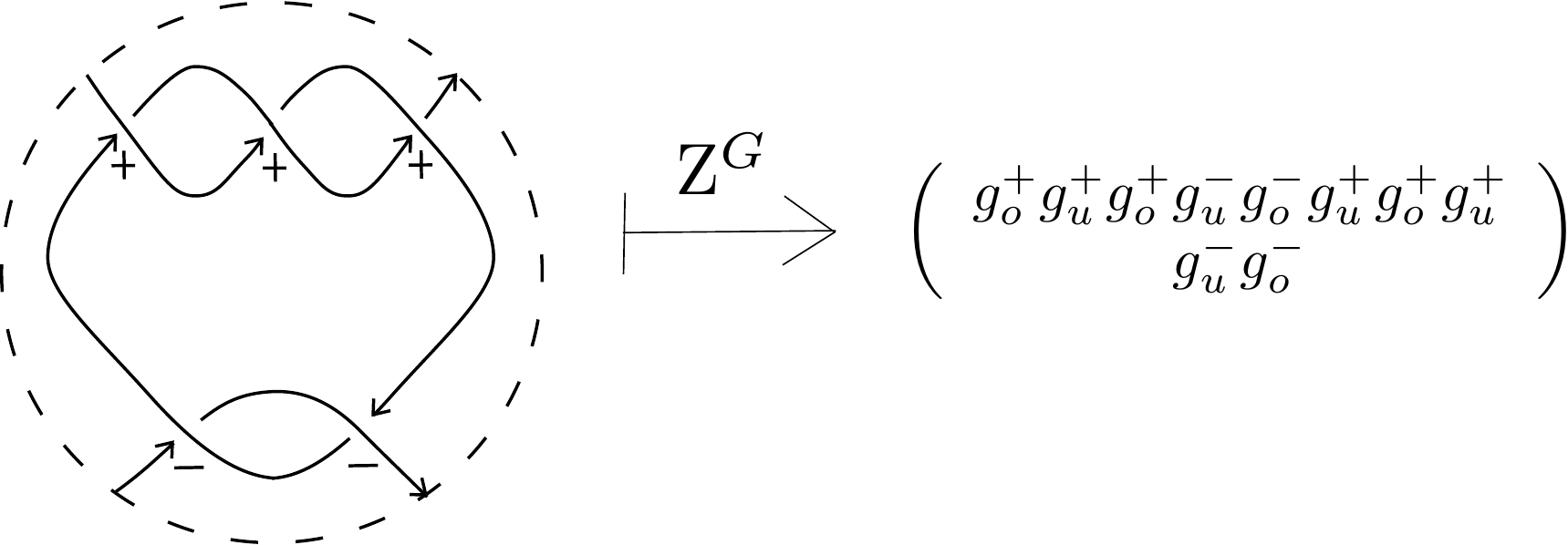}
    \caption{Computing $Z^G$ of a tangle: (a) assigning values to crossings. (b) collecting along strands}\label{fig:MonoidZ}
  \end{figure}
\else
  \begin{figure}[h]
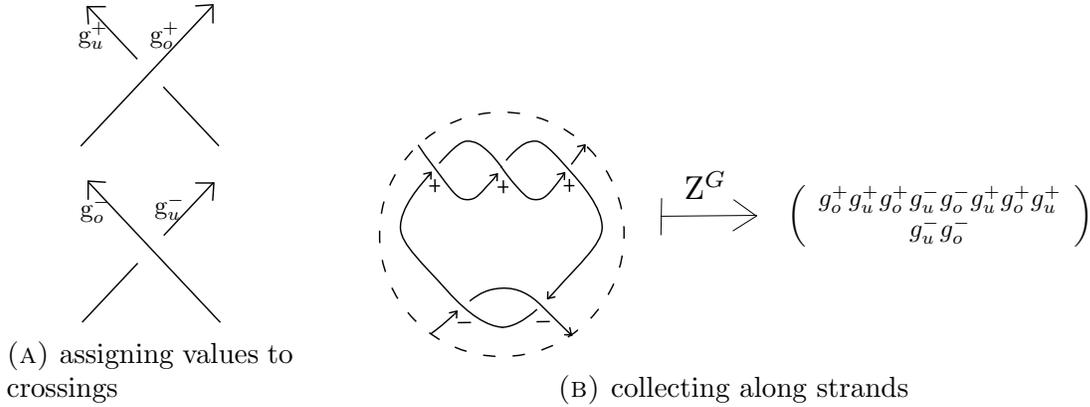

    \centering
    \begin{subfigure}[b]{1.5in}
      \centering
      \includegraphics[scale=0.54]{Crossings.pdf}
      \caption{assigning values to crossings}
    \end{subfigure}
    \begin{subfigure}[b]{4.5in}
      \centering
      \includegraphics[scale=0.54]{Zfigure.pdf}
      \caption{collecting along strands}
    \end{subfigure}
    \caption{Computing $Z^G$ of a tangle}\label{fig:MonoidZ}
  \end{figure}
\fi

Unfortunately, the gods are not so kind and $Z^G$ is not worth much more than the effort that went in it.  Indeed,
invariance under the Reidemeister $I\!I$ move (see Figure 2) demands $g_o^-={(g_o^+)}^{-1}$ and $g_u^-={(g_u^+)}^{-1}$, while Reidemeister $I\!I\!I$ adds that $g_o^+$ and $g_u^+$, as well as $g_o^-$ and $g_u^-$, commute. As a result, every component of $Z^G(T)$ collapses to the form $g_o^ag_u^b$ for some integers $a$ and $b$, so all the information to bring home is the signed number of times a given strand crosses over or under other strands.
It will turn out, nevertheless, that a generalized version of this procedure yields an amply non-trivial invariant with novel properties.

\if\jktr y
  \begin{figure}[htpb]\centering
    \includegraphics[scale=0.38]{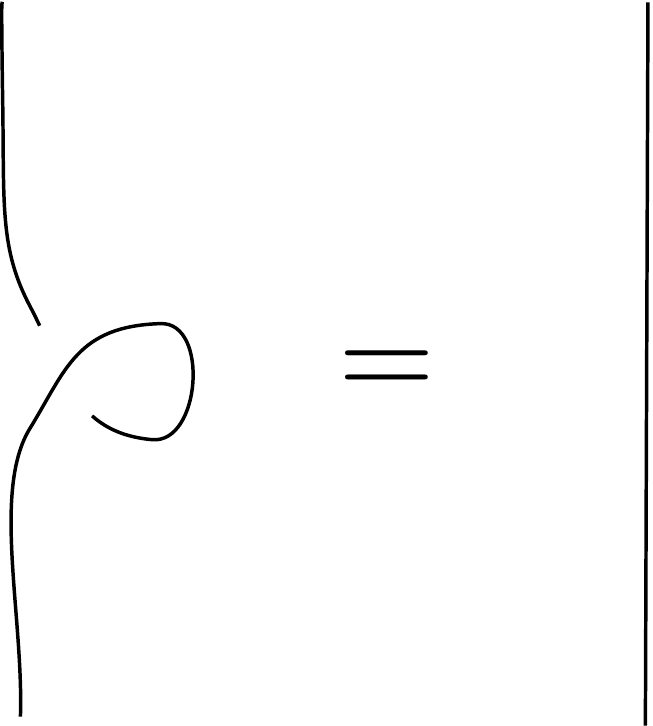}
    \qquad
    \includegraphics[scale=0.38]{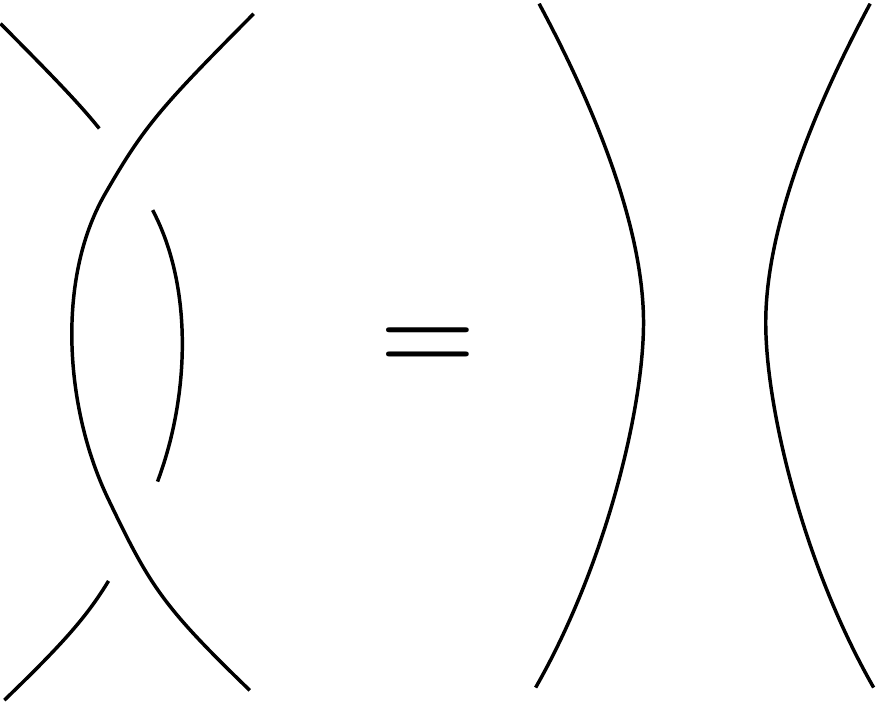}
    \qquad
    \includegraphics[scale=0.38]{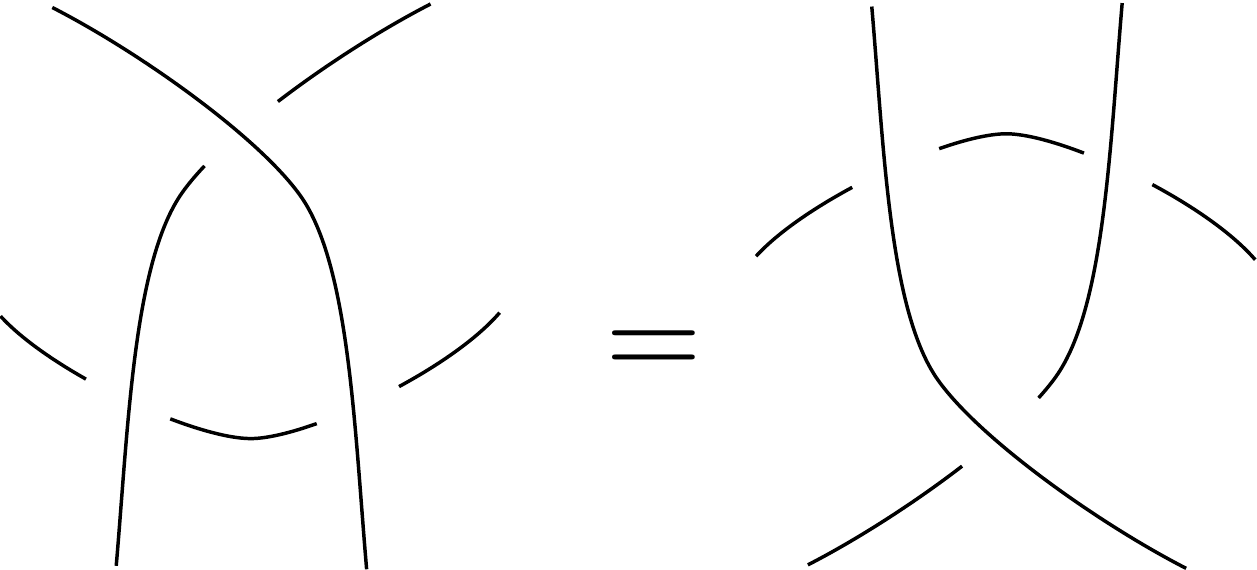}
    \caption{The three Reidemeister moves: $I$, $I\!I$, $I\!I\!I$}
  \end{figure}
\else
  \begin{figure}[htpb]
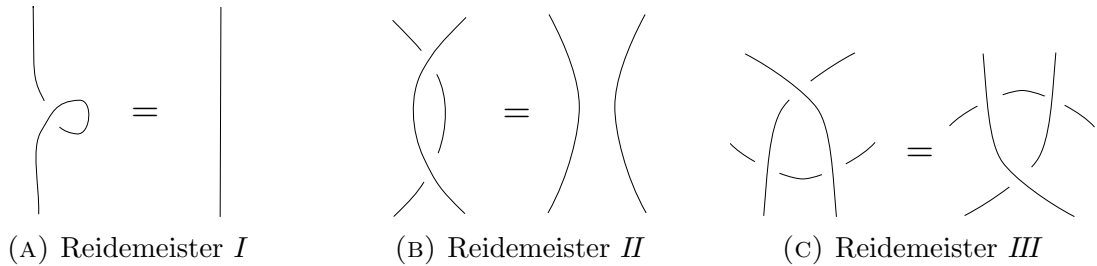

    \centering
    \begin{subfigure}[b]{2in}
      \centering
      \includegraphics[scale=0.38]{R1.pdf}
      \caption{Reidemeister $I$}
    \end{subfigure}
    \begin{subfigure}[b]{2in}
      \centering
      \includegraphics[scale=0.38]{R2.pdf}
      \caption{Reidemeister $I\!I$}
    \end{subfigure}
    \begin{subfigure}[b]{2in}
      \centering
      \includegraphics[scale=0.38]{R3.pdf}
      \caption{Reidemeister $I\!I\!I$}
    \end{subfigure}
    \caption{The three Reidemeister moves}
  \end{figure}
\fi

\section{A better invariant: $Z^\beta$}\label{sec:Zbeta}

The invariant that we wish to introduce can be thought of as taking values in a meta-monoid. This is a generalization of what we call a ``monoid computer'':

\subsection{Preliminary: A Monoid Computer}\label{ssec:preliminary}

If $X$ is a finite set and $G$ is a monoid we let $G^X$ denote the set of all possible assignments of elements of $G$ to the set X; these are \textquoteleft\textquoteleft$G$-valued datasets, with registers labelled by the elements of $X$\textquoteright\textquoteright.
\begin{figure}[h]
\centering
\includegraphics[scale=0.5]{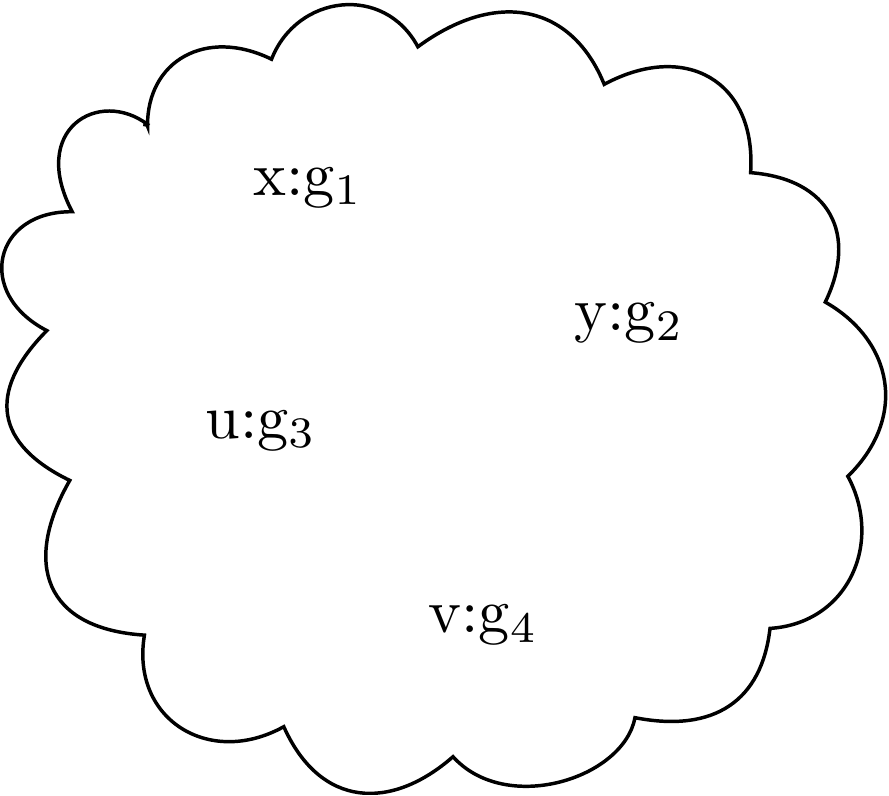}
\caption{A typical element of $G^{\{x,y,u,v\}}$}
\end{figure}

A monoid computer can manipulate registers in some prescribed ways.
For example, if $X$ does not contain $x$, $y$ and $z$, define \space $m_z^{xy}\colon G^{X\cup \{x,y\}} \to G^{X\cup \{z\}}$ using the monoid multiplication, $\{x\colon g_1, y\colon g_2\} \mapsto \{z\colon g_1g_2\}$.
There are obvious operations for renaming or deleting a register, and inserting the identity in a new register, respectively denoted $\rho_y^x$, $d^x$ and $e_y$, and respectively implemented on $G^{X\cup \{x\}}$ by fixing the content of $X$ and mapping $\{x\colon g\}$ to $\{y\colon g\}$, $\{\}$ and $\{x\colon g,y\colon e\}$.
In addition there is a binary operation for merging data sets, $\bigcup\colon G^X \times G^Y \to G^{X \cup Y}$, which takes two data sets $P$ and $Q$ and forms their disjoint union $P \cup Q$. We can compose the aforementioned maps if labels match correctly, and we do so from left to right with the aid of the notation $\sslash$. For example, we write $P \sslash \rho^x_y\sslash \rho^y_z$ to rename the register $x$ of P first to $y$, then to $z$.

\subsection{Meta-Monoids}\label{ssec:metamonoids}
The operations on a monoid computer obey a certain set of basic set-theoretic axioms as well as axioms inherited from the monoid $G$.
A meta-monoid is an abstract computer that satisfies some but not all of those axioms. We postpone the precise definition to Section~\ref{sec:metamonoids}. It may be best to begin with examples and a prototypical one is as follows. Let $G_X:=M_{X\times X}(\textbf{Z})$ denote (not in reference to any monoid $G$) the set of $|X| \times |X|$ matrices of integers with rows and columns labelled by $X$. The operation of \textquoteleft\textquoteleft multiplication\textquoteright\textquoteright, on say, $3\times3$ matrices, $m_z^{xy}\colon G_{\{x,y,w\}}\to G_{\{z,w\}}$, is defined by simultaneously adding rows and columns labelled by $x$ and $y$:
\begin{equation*}
 \bordermatrix{
& x &
y &
w \cr
x& a & b & c \cr
y& d & e & f \cr
w& g & h & i
 }
 \mapsto
 \bordermatrix{
 &z&w\cr
 z & a + b + d + e & c + f \cr
 w& g + h & i}
\end{equation*}

While still satisfying the associativity condition $m^{xy}_u\sslash m^{uv}_w = m^{yv}_u\sslash m^{xu}_w$, this example differs from a monoid computer by the failure of a critical axiom: if $P \in G_{\{x,y\}}$, $$d_yP\cup d_xP \ne P$$ Indeed, if $P \in G_{\{x,y\}}$ is the matrix $\bordermatrix{&x&y\cr x&a&b\cr y&c&d}$,  then $$d_yP\cup d_xP=\bordermatrix{&x&y\cr x&a&0\cr y&0&d}\ne P$$

\subsection{Meta-Bicrossed Products}

Suppose a group $G$ is given as the product $G=TH$ of two of its subgroups, where $T\cap H=\{e\}$. Then also $G=HT$ \footnotemark \footnotetext{Indeed, if $g^{-1}=th$, then $g=h^{-1}t^{-1}$, so $g^{-1}\in TH$ implies $g\in HT$, and as $TH=G$, also $HT=G$.} and every element of $G$ has unique\footnotemark  \space \footnotetext{Separation of variables: suppose $g=h_1t_1=h_2 t_2$. Then we have $h_2^{-1}h_1=t_2t_1^{-1}$, which implies that $h_1=h_2$ and $t_1=t_2$ since $h_2^{-1}h_1\in H$, $t_2t_1^{-1}\in T$, and $H\cap T={\{e\}}$.}representations of the form $th$ and $h^\prime t^\prime$ where $h,h^\prime\in H$ and $t, t^\prime\in T$.  Accordingly there is a \textquoteleft\textquoteleft  swap\textquoteright\textquoteright\space  map $sw\colon T\times H\to H \times T$, $(t,h)\mapsto(h^\prime,t^\prime)$ such that if $g=th$ then $g=h^\prime t^\prime$ also. The swap map satisfies some relations; in monoid-computer language, the important ones are as in Figure~\ref{fig:swap}. Conversely, provided that the swap map satisfies the relations in Figure~\ref{fig:swap}, the data $(H,T,sw)$ determines a monoid $G$, with product given by $\{(h_1,t_1), (h_2, t_2)\} \mapsto (h_1 h_2^\prime, t_1^\prime t_2)$ where $sw(t_1,h_2)=(h_2^\prime, t_1^\prime)$. $G$ is called the bicrossed product of $H$ and $T$, which we could denote $(H\times T)_{sw}$. In a semidirect product, one of $H$ or $T$ is normal (say $T$) and the swap map is $sw\colon(t,h) \mapsto (h, h^{-1} t h)$.

\if\jktr y
  \begin{figure}\centering
    \includegraphics[scale=0.5]{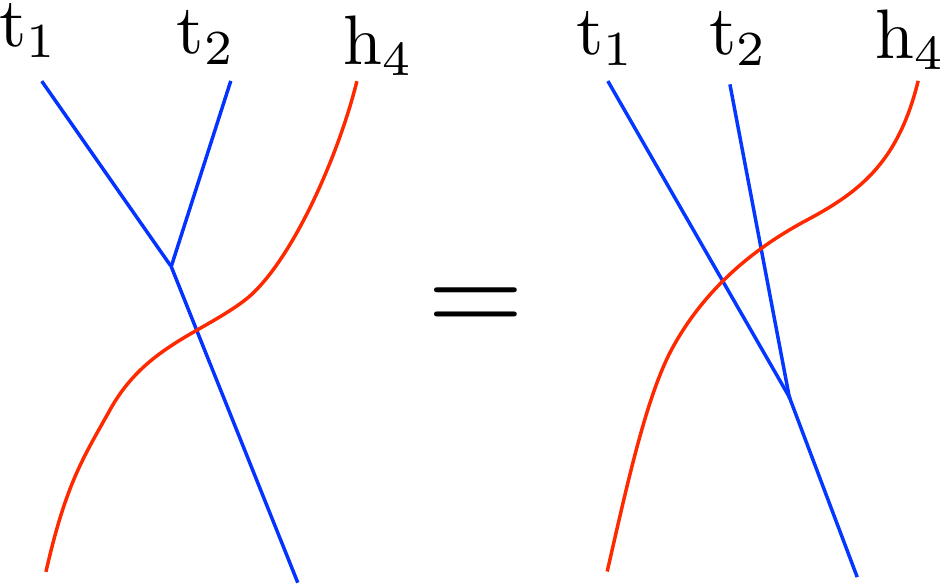}
    \qquad
    \includegraphics[scale=0.5]{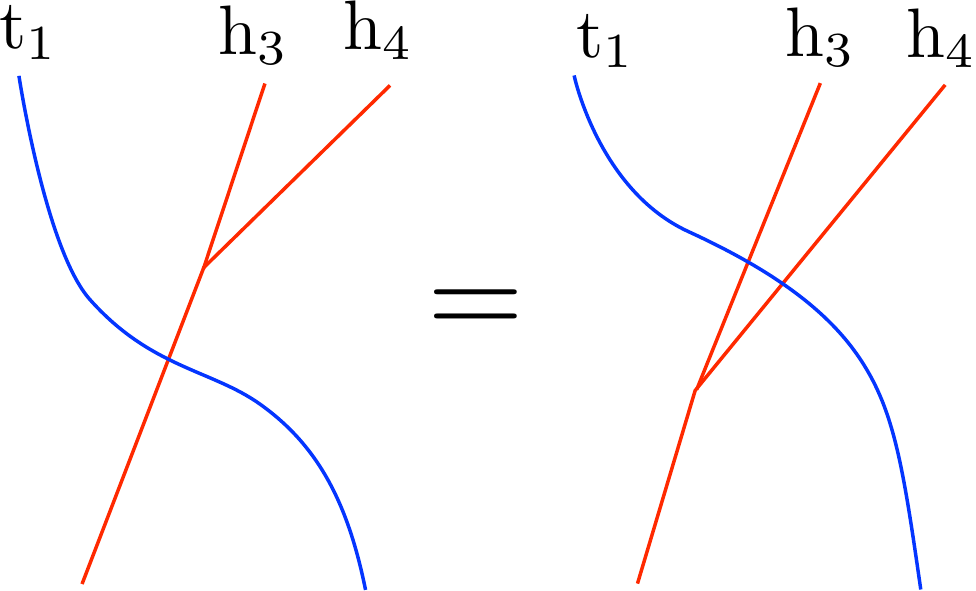}
    \caption{swap operation axioms. $tm$ and $hm$ stand for multiplication in $T$ and $H$ respectively: (a) $tm_1^{12} \sslash  sw_{14}=sw_{24}\sslash sw_{14}\sslash tm_1^{12}$. (b) $hm_3^{34}\sslash sw_{13}=sw_{13}\sslash sw_{14}\sslash hm_3^{34}$} \label{fig:swap}
  \end{figure}
\else
  \begin{figure}
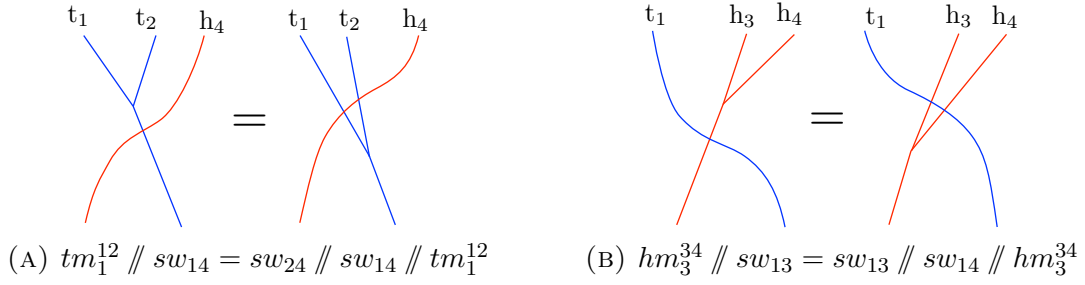

    \centering
    \begin{subfigure}[b]{3in}
      \centering
      \includegraphics[scale=0.5]{sw1.pdf}
      \caption{$tm_1^{12} \sslash  sw_{14}=sw_{24}\sslash sw_{14}\sslash tm_1^{12}$} \label{eq:sw1}
    \end{subfigure}
    \begin{subfigure}[b]{3in}
      \centering
      \includegraphics[scale=0.5]{sw2.pdf}
      \caption{$hm_3^{34}\sslash sw_{13}=sw_{13}\sslash sw_{14}\sslash hm_3^{34}$} \label{eq:sw2}
    \end{subfigure}
    \caption{swap operation axioms. $tm$ and $hm$ stand for multiplication in $T$ and $H$ respectively.} \label{fig:swap}
  \end{figure}
\fi

The corresponding notion of a meta-bicrossed product is a collection of sets $\beta(\eta,\tau)$ indexed by all \textit{pairs} of finite sets $\eta$ and $\tau$ ($\eta$ for \textquoteleft\textquoteleft heads\textquoteright\textquoteright, $\tau$ for \textquoteleft\textquoteleft tails\textquoteright\textquoteright), and equipped with multiplication maps $tm_z^{xy}$ ($x$, $y$ and $z$ tail labels), $hm_z^{xy}$ ($x$, $y$ and $z$ head labels), and a swap map $sw_{xy}^{th}$ (where $t$ and $h$ indicate that $x$ is a tail label and $y$ is a head label --- note that $sw_{yx}^{ht}$ is in general a different map) satisfying (a) and (b).

\parpic[r]{$
\xymatrix@R=0.25in{
  & \parbox{1.2in}{\centering$\bordermatrix{
  & h_1 \cr
  t_1 & a + b \cr
  t_2 & c+d \cr
  t_3 & e+f
}$}
    \\
\parbox{1in}{\centering$\bordermatrix{
& h_1 & h_2  \cr
t_1& a & b  \cr
t_2& c & d\cr
t_3& e & f
 }$}
\ar[ru]^{hm_1^{1,2}} \ar[r]^{tm_1^{1,2}} \ar[rd]_{sw_{1,2}^{th}}
&
 \parbox{1.4in}{\centering$\bordermatrix{
 &h_1 & h_2\cr
 t_1 & a+c & b+d \cr
 t_3 & e & f
}$}
\\
   & \parbox{1.2in}{$\centering\bordermatrix{
& h_1 & h_2  \cr
t_1 & a & b  \cr
t_2 & c & d \cr
t_3 & e & f
 }$}
}
$}
Given the above we can make a ``monoid multiplication'' map out of the head and tail multiplication maps via $gm_z^{xy}:= sw_{xy}^{th} \sslash tm_z^{xy} \sslash hm_z^{xy}$. Thus a meta-bicrossed product defines a meta-monoid with $\Gamma_X=\beta(X,X)$. An example of a meta-bicrossed product is given by the rectangular matrices, $\mu(\eta,\tau) := M_{\tau\times\eta}(\textbf{Z})$, with $tm_z^{xy}$ and $hm_z^{xy}$ corresponding to adding two rows and adding two columns, and swap being the trivial operation. Here $\Gamma_X$ is the same as the first example of Section~\ref{ssec:metamonoids}. An example with a non-trivial swap map will shortly follow.

\subsection{$\beta$ Calculus}\label{ssec:betacalculus}

The $\beta$ calculus has an arcane origin \cite{wclips}\footnote{in which, among other things, the ``heads and tails'' vocabulary is motivated.} which we will not discuss. We expect that it can be presented in a much simpler and fitting context than that in which it was discovered. Accordingly we will simply pull it out of a hat. Though note that many of our formulas bear close resemblance to formulas in~\cite{LeDimet:Gassner, KirkLivingstonWang:Gassner, CimasoniTuraev:LagrangianRepresentation}.

Let $\beta(\eta,\tau)$ be (again, in reference to \textit{sets} $\eta$ and $\tau$) the collection of arrays with rows labeled by $t_i \in \tau$ and columns labeled by $h_j \in \eta$, along with a distinguished element $\omega$. Such arrays are conveniently presented in the following format:

\begin{equation*}
\begin{array}{c | c c c}
\omega  & h_1 & h_2 & \dots \\
\hline
 t_1 & \alpha_{11} & \alpha_{12} & \cdot \\
 t_2 & \alpha_{21} & \alpha_{22} & \cdot \\
 \vdots & \cdot & \cdot & \cdot \\
\end{array}
\end{equation*}

The $\alpha_{ij}$ and $\omega$ are rational functions of variables $T_i$, which are in bijection with the row labels $t_i$.

$\beta(\eta,\tau)$ is equipped with a peculiar set of operations. Despite being repulsive at sight, they are completely elementary. They are defined as follows:

\vspace{10pt}

\noindent\begin{minipage}{0.58\linewidth}
\begin{multline*}
\shoveright{
tm_z^{xy}\colon
\begin{array}{c|c}
 \omega  & \dots \\
\hline
 t_x & \alpha  \\
 t_y & \beta  \\
\vdots & \gamma  \\
\end{array}
\mapsto
\begin{array}{c|c}
 \omega  & \dots \\
\hline
 t_z & \alpha +\beta  \\
 \vdots & \gamma  \\
\end{array}
}
\end{multline*}
\end{minipage}
\begin{minipage}{0.38\linewidth}
Here $\alpha$ and $\beta$ are rows and $\gamma$ is a matrix. The sum $\alpha+\beta$ is accompanied by the corresponding change of variables $T_x$, $T_y$ $\mapsto T_z$.
\end{minipage}

\vspace{5pt}

\noindent\begin{minipage}{0.58\linewidth}
\begin{multline*}
\shoveright{
hm_z^{xy}\colon
\begin{array}{c|ccc}
 \omega  & h_x & h_y & \dots \\
\hline
 \vdots & \alpha  & \beta  & \gamma  \\
\end{array}
\mapsto
\begin{array}{c|cc}
 \omega  & h_z & \dots \\
\hline
 \vdots& \alpha+\beta+\langle \alpha \rangle \beta & \gamma  \\
\end{array}
}
\end{multline*}
\end{minipage}\hfill
\begin{minipage}{0.33\linewidth}
Here $\alpha$ and $\beta$ are columns, $\gamma$ is a matrix, and $\langle \alpha \rangle =  \sum_i \alpha_i$.
\end{minipage}

\vspace{5pt}

\noindent\begin{minipage}{0.58\linewidth}
\begin{multline*}
\shoveright{
sw_{xy}^{th}\colon
\begin{array}{c|cc}
 \omega  & h_y & \dots \\
\hline
 t_x & \alpha  & \beta  \\
 \vdots & \gamma  & \delta  \\
\end{array}
\mapsto
\begin{array}{c|cc}
 \omega \epsilon  & h_y & \dots \\
\hline
 t_x &\alpha (1+\langle \gamma \rangle / \epsilon ) & \beta (1+\langle \gamma \rangle /\epsilon ) \\
 \vdots & \gamma /\epsilon & \delta -\gamma \beta /\epsilon \\
\end{array}
}
\end{multline*}
\end{minipage}\hfill
\begin{minipage}{0.33\linewidth}
Here $\alpha$ is a single entry, $\beta$ is a row, $\gamma$ is a column, and $\delta$ is a matrix comprised of the rest. $\epsilon=1+\alpha$. Note also that $\gamma \beta$ is the matrix product of the column $\gamma$ with the row $\beta$ and hence has the same dimensions as the matrix $\delta$.
\end{minipage}

\vspace{20pt}

We also need the disjoint union, defined by

\begin{equation*}
\begin{array}{c|c}
\omega_1 & H_1 \\
\hline
T_1 & \alpha_1 \\
\end{array}
\cup
\begin{array}{c|c}
\omega_1 & H_1 \\
\hline
T_1 & \alpha_1 \\
\end{array}
=
\begin{array}{c|cc}
\omega_1\omega_2 & H_1 & H_2 \\
\hline
T_1 & \alpha_1 & 0 \\
T_2 & 0 & \alpha_2 \\
\end{array}
\end{equation*}

We make $\beta$ into a meta-monoid via the ``monoid-multiplication'' map $gm^{xy}_z:=sw^{th}_{xy}\sslash tm_z^{xy}\sslash hm_z^{xy}$. We will later set out to make proper definitions, write down the remaining operations, and establish the following

\begin{theorem}
$\beta$ is a meta-bicrossed product.
\end{theorem}

Finally there are two elements which will serve as a pair of ``R-matrices'', analogous to the pair of pairs $(g_o^\pm,g_u^\pm)$ of $Z^G$:

\begin{align*}
R_{xy}^+ =
\begin{array}{c|cc}
1 & h_x & h_y \\
\hline
t_x & 0 & T_x-1 \\
t_y & 0 & 0 \\
\end{array} & &
R_{xy}^- =
\begin{array}{c|cc}
1 & h_x & h_y \\
\hline
t_x & 0 & T_x^{-1}-1 \\
t_y & 0 & 0 \\
\end{array}
\end{align*}


\subsection{$Z^\beta$}
Let $T$ be again an oriented tangle diagram.
At each crossing, assign a number to the upper strand and to the lower strand. Using the $R^\pm_{xy}$ of above, form the disjoint union $\bigcup_{\{i,j\}} R_{ij}^\pm$ where $\{i,j\}$ runs over all pairs assigned to crossings, with $i$ labelling the upper strand and $j$ labelling the lower strand, and where $\pm$ is determined by the sign of the given crossing. Now for each strand multiply all the labels in the order in which they appear. That is, if the first label on the strand is $k$, repeatedly apply $gm^{kl}_k$ where $l$ runs over all labels subsequently encountered on the strand (in order). If $T$ has $n$ strands, the result is an $n\times n$ array with an extra corner element. Call this array $Z^\beta(T)$. Those were a lot of words, so take for example the knot $8_{17}$ illustrated in Figure \ref{fig:817}. In this case, form the disjoint union\footnote{From now on we omit the $\cup$ in disjoint unions: $\beta_1\beta_2:=\beta_1\cup\beta_2$.}

\begin{equation*}
R_{12,1}^-R_{2,7}^-R_{8,3}^-R_{4,11}^-R_{16,5}^+R_{6,13}^+R_{14,3}^+R_{10,15}^+,
\end{equation*}

which is given by the following array\footnote{We suppress rows/columns of zeros.}:

\begin{equation*}
\begin{array}{c|cccccccc}
1 & h_1 & h_3 & h_5 & h_7 & h_9 & h_{11} & h_{13} & h_{15} \\
\hline
t_2 & 0 & 0 & 0 & T_2^{-1} - 1 & 0 & 0 & 0 & 0 \\
t_4 & 0 & 0 & 0 & 0 & 0 & T_4^{-1} - 1 & 0 & 0 \\
t_6 & 0 & 0 & 0 & 0 & 0 & 0 & T_6 - 1 & 0 \\
t_8 & 0 & T_ 8^{-1} - 1 & 0 & 0 & 0 & 0 & 0 & 0 \\
t_{10} & 0 & 0 & 0 & 0 & 0 & 0 & 0 & T_{10} - 1\\
t_{12} & T_{12}^{-1} - 1 & 0 & 0 & 0 & 0 & 0 & 0 & 0\\
t_{14} & 0 & 0 & 0 & 0 &T_{14}-1 & 0 & 0 & 0 \\
t_{16} & 0 & 0 & T_{16} - 1 & 0 & 0 & 0 & 0 & 0 \\
\end{array}
\end{equation*}

Then apply the multiplications $gm^{1k}_1$, with $k$ running from 2 to 16, to get the following $1\times 1$ array with corner element:
\begin{equation*}
\begin{array}{c|c}  -T_1^{-3}+4T_1^{-2}-8T_1^{-1}+11-8T_1+4T_1^2-T_1^3 & h_1 \\ \hline t_1 & 0 \end{array}
\end{equation*}


\begin{theorem}
$Z^\beta$ is an invariant of oriented tangle diagrams.
\end{theorem}
\begin{proof}
Straightforward check. We do the computation for the Reidemeister $I\!I\!I$ move to illustrate.
The disjoint unions for each side of the equality are given by:
\begin{center}
\begin{minipage}{0.25\linewidth}
\centering
\includegraphics[scale=0.4]{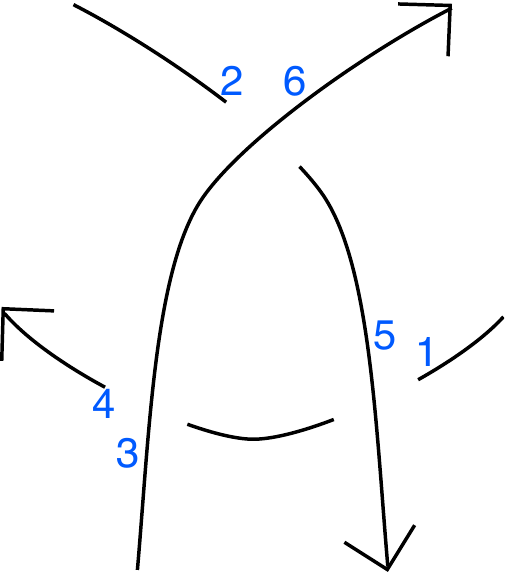}

\vspace{20pt}

\includegraphics[scale=0.4]{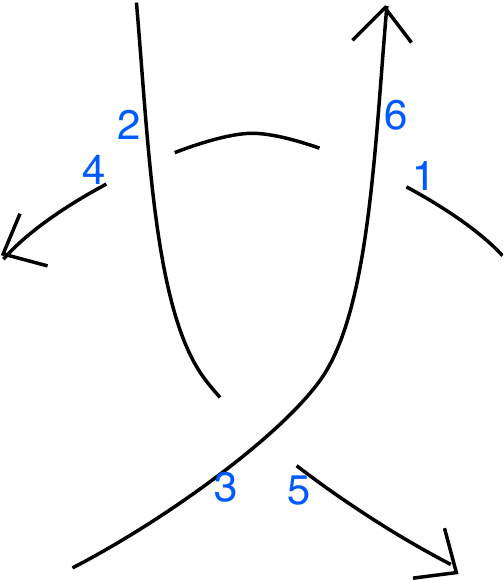}
\end{minipage}
\begin{minipage}{0.5\linewidth}
\begin{equation*}
R^-_{1,5}R^-_{6,2}R^+_{3,4}=
\begin{array}{c|ccc}
1 & h_1 & h_2 & h_4 \\
\hline
t_3 & 0 & 0 & T_3 - 1 \\
t_5 & T_5^{-1} - 1 & 0 & 0 \\
t_6 & 0 & T_6^{-1} - 1 & 0
\end{array}
\end{equation*}
\vspace{20pt}
\begin{equation*}
R^+_{6,1}R^-_{2,4}R^-_{3,5}=
\begin{array}{c|ccc}
1 & h_1 & h_4 & h_5 \\
\hline
t_3 & 0 & T_2^{-1} -1 & 0 \\
t_5 & 0 & 0 & T_3^{-1} -1 \\
t_6 & T_6 -1 & 0 & 0
\end{array}
\end{equation*}
\end{minipage}
\end{center}

Then one checks that indeed
\begin{multline*}
 R^-_{1,5}R^-_{6,2}R^+_{3,4}\sslash gm_1^{1,4} \sslash gm_2^{2,5} \sslash gm_3^{3,6}=
 R^+_{6,1}R^-_{2,4}R^-_{3,5}\sslash gm_1^{1,4} \sslash gm_2^{2,5} \sslash gm_3^{3,6} \\
 =
 \begin{array}{c|cc}
 1 & h_1 & h_2 \\
 \hline
 t_1 & T_2^{-1} - 1& 0 \\
 t_2 & T_2^{-1}(T_3-1) & T_3^{-1}-1
 \end{array}
\end{multline*}
\end{proof}

One philosophically appealing major property of $Z^\beta$ is that the operations used to compute it have a literal interpretation of gluing crossings together. In particular, at every stage of the computation we get an invariant of the tangle\footnotemark \space made of all the crossings but only those for which the corresponding $gm$ was carried out have been glued. Additionally, unlike other existing extensions of the Alexander polynomial to tangles, $Z^\beta$ takes values in spaces of polynomial size, at every step of the calculation.
\footnotetext{The careful reader may wish to peek ahead at Section \ref{ssec:vtangles} for a better grasp of this statement.}

\subsection{Knots and links}\label{ssec:links}

\begin{conjecture} \label{conj:Alexander}
Restricted to long knots (which are the same as round knots), the corner element of $Z^\beta$ is the Alexander polynomial. Restricted to string links (which map surjectively to links), $Z^\beta$ contains the multivariable Alexander polynomial.
\end{conjecture}

While we are shy of a formal proof, the computer evidence behind Conjecture~\ref{conj:Alexander} is overwhelming. See Section~\ref{ssec:KnotsAndLinks}.

\section{More on meta-monoids}\label{sec:metamonoids}


\subsection{The meta-monoid of coloured v-tangles}\label{ssec:vtangles}
When one tries to follow the interpretation of the computation of $Z^\beta$ as progressively attaching crossings together to form a tangle, one will in general encounter a step where the tangle becomes non-planar (a strand will have to go through another in an ``artificial'' crossing to reach the boundary disk). See Figure~\ref{fig:817}. Such tangles are called virtual or v-tangles and constitute a rich subject of study on their own; see~\cite{Kauffman:VirtualKnotTheory}. We will be content with acknowledging their existence and giving them a name.

\if\jktr y
  \begin{figure}[h]
    \includegraphics[scale=0.45]{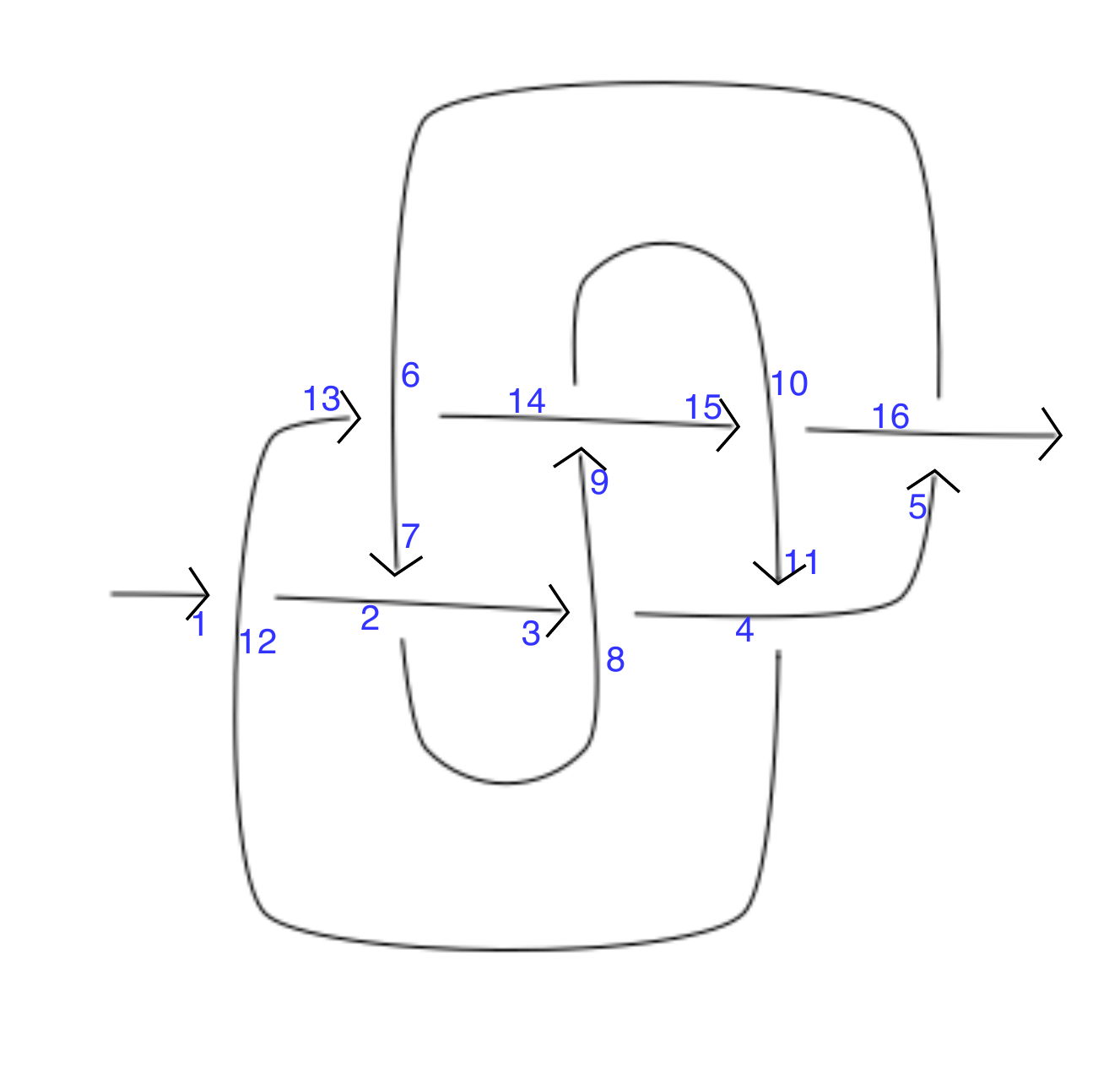}
    \quad
    \includegraphics[scale=0.45]{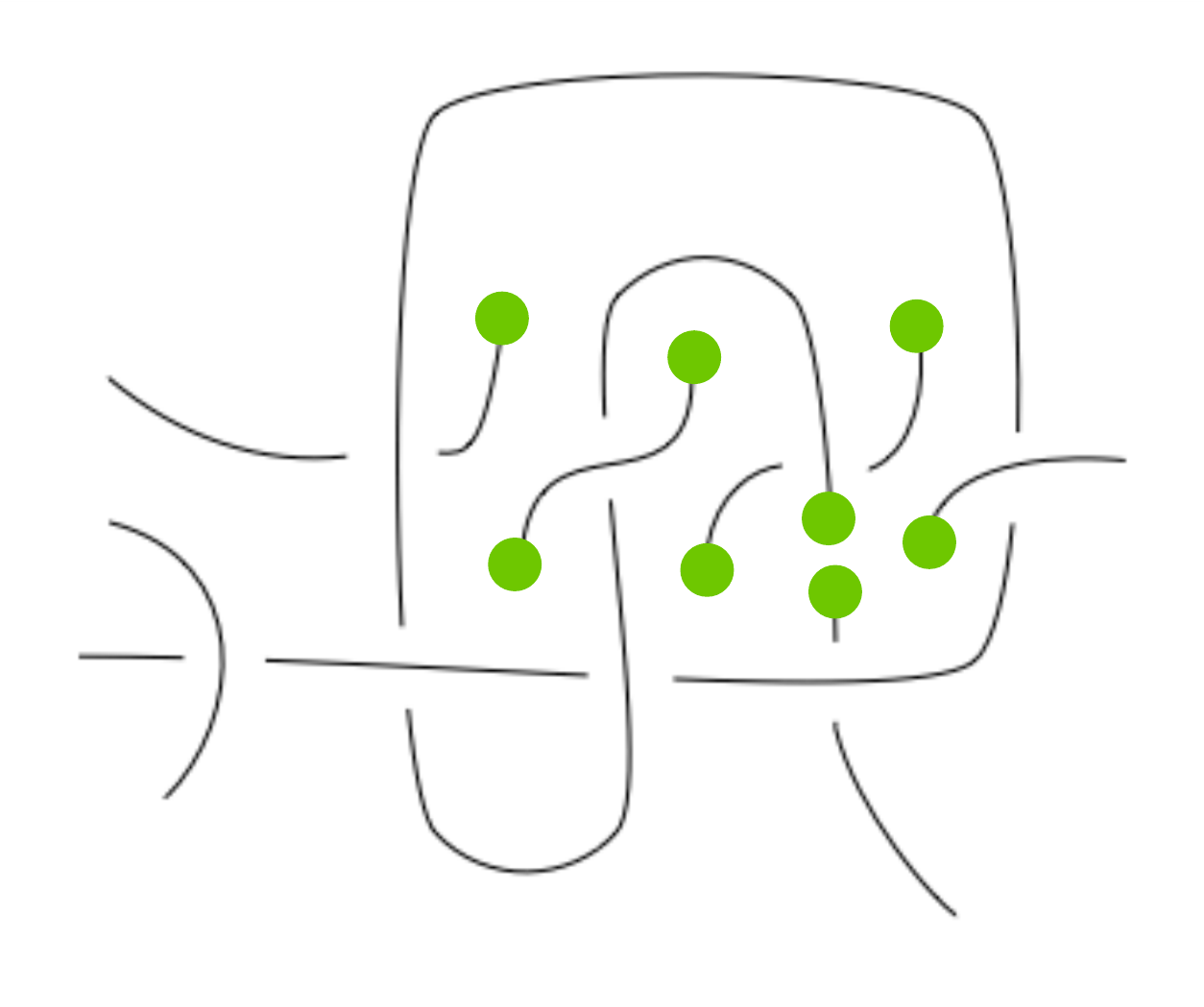}
    \caption{The knot $8_{17}$: (a) With crossings labelled. (b) After attaching crossings 1 through 10. The arcs with green dots can not make it out to the boundary disk.} \label{fig:817}
  \end{figure}
\else
  \begin{figure}[h]
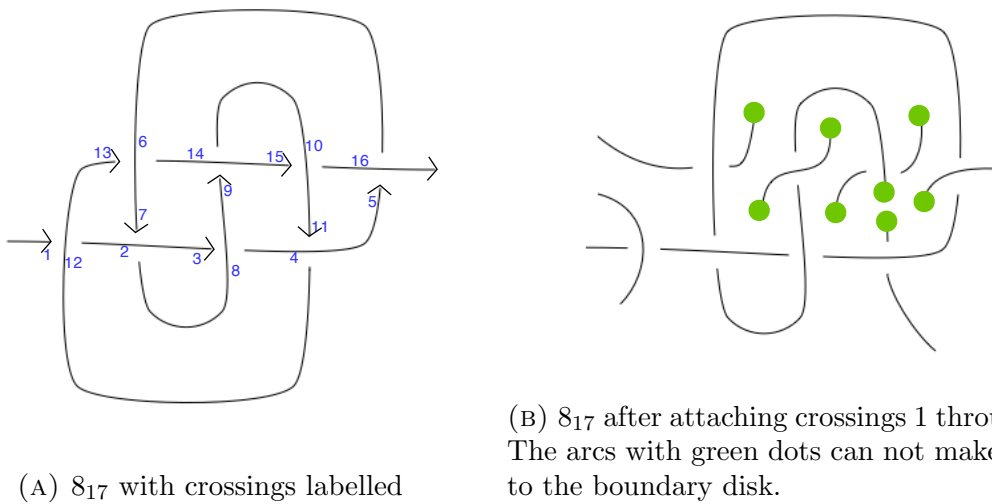

    \centering
    \begin{subfigure}[b]{3in}
      \centering
      \includegraphics[scale=0.5]{817labeled.pdf}
      \caption{$8_{17}$ with crossings labelled}
    \end{subfigure}
    \begin{subfigure}[b]{3in}
      \centering
      \includegraphics[scale=0.5]{817_10x_greened.pdf}
      \caption{$8_{17}$ after attaching crossings 1 through 10. The arcs with green dots can not make it out to the boundary disk.}
    \end{subfigure}
    \caption{The knot $8_{17}$} \label{fig:817}
  \end{figure}
\fi

If $X$ is a finite set, oriented $X$-coloured pure\footnotemark \footnotetext{Pure means that the tangles have no closed component.} virtual tangles form a meta-monoid. The operation $m^{xy}_z$ attaches the head of strand $x$ to the tail of strand $y$ (possibly through a few virtual crossings) and names the resulting strand $z$\footnote{Remark: this is \textit{not} a meta-generalization of the group structure on braids.}.

\subsection{Some familiar invariants}

We have already suggested that $Z^G$ and $Z^{\beta}$ take values in meta-monoids. Some more traditional invariants can also be cast in meta-monoid context. Note that $Z^G$ is in fact very traditional, being nothing more than linking numbers. We invite the reader familiar with the fundamental group of the complement of a tangle to consider the following set-up:

Let $G_{\{x_1, \dots, x_n\}}=\{ ( \Gamma, m_1, l_1, \dots , m_n, l_n); \text{$\Gamma$ is a group}; m_i, l_i \in \Gamma \}$. The multiplication map that corresponds to what happens to the meridians and longitudes when one plugs a strand into another is
\[
  m^{ij}_i(\Gamma, m_1,l_1, \dots, m_n,l_n)=(\Gamma / (m_j=l_i^{-1}m_i l_i), m_1, l_1 l_2, \dots, \widehat{m_j}, \widehat{l_j}, \dots,  m_n, l_n)
\]

Also the fundamental group of the complement of two disjoint tangles is the free product of the respective fundamental groups, so we define also
\begin{multline*}
  (\Gamma^1, m_1^1,l_1^1, \dots, m_n^1,l_n^1) \cup (\Gamma^2, m_1^2,l_1^2, \dots, m_k^2,l_k^2) \\
  =(\Gamma^1\star \Gamma^2, m_1^1,l_1^1,\dots, m_n^1,l_n^1,  m_1^2,l_1^2,\dots, m_k^2,l_k^2).
\end{multline*}

\subsection{Definitions}\label{ssec:def}

We now proceed to laying down the details of the definitions of meta-monoids and meta-bicrossed products.

A meta-monoid is a collection of sets $\Gamma$ indexed by all finite sets, equipped with operations $m^{xy}_z\colon \Gamma_{\{x,y\}\cup X} \to \Gamma_{\{z\}\cup X}$, $e_x\colon \Gamma_X \to \Gamma_{\{x\} \cup X}$, $d_x\colon \Gamma_{\{x\} \cup X} \to \Gamma_X$, and $\bigcup\colon \Gamma_X \times \Gamma_Y \to \Gamma_{X \cup Y}$ satisfying the following:

\vskip 3mm

\noindent\begin{minipage}[t]{0.45\linewidth}
``Monoid theory'' axioms
\begin{itemize}

\item{$e_x\sslash m_z^{xy}=\rho_z^y$} (left identity)
\item{$e_y\sslash m_z^{xy}=\rho_z^x$} (right identity)
\item{$m_u^{xy}\sslash m_v^{uz}=m_u^{yz}\sslash m_v^{xu}$} (associativity)
\end{itemize}
\end{minipage}\hfill
\begin{minipage}[t]{0.5\linewidth}
``Set manipulation'' axioms
\begin{itemize}

\item{$\rho_x^y \sslash \rho_y^x=id$}
\item{$\rho_y^x\ \sslash \rho_z^y=\rho_z^x$}
\item{$\rho_y^x \sslash d_y=d_x$}
\item{$m^{xy}_z \sslash d_z = d_x \sslash d_y$}
\item{$e_x \sslash d_x = id$}
\item{$m^{xy}_z \sslash \rho_u^z = m^{xy}_u$}
\item{$\rho^x_u \sslash m^{uy}_z = m^{xy}_z$}
\item{$e_x \sslash \rho_y^x = e_y$}
\item Operations involving disjoint sets of labels commute (e.g. $e_x \sslash e_y = e_y\sslash e_x$)
\end{itemize}
\end{minipage}

\vskip 3mm

A meta-bicrossed product is a collection of sets $\Gamma$ indexed by all pairs of finite sets, equipped with maps $hm$, $tm$, and $sw$, such that:
\begin{itemize}
\item{$hm^{xy}_{z}\colon\Gamma(\eta\cup \{x,y\},\tau_0) \to \Gamma(\eta\cup \{z\}, \tau_0)$ and $tm^{x y}_{z}\colon\Gamma(\eta_0,\tau\cup \{x,y\}) \to \Gamma(\eta_0, \tau\cup \{z\})$ define a meta-monoid structure for each fixed choice of $\tau_0$ and $\eta_0$, respectively.}
\item{$sw_{xy}$ satisfies the following relations (recall Figure \ref{fig:swap})}
\begin{itemize}
\item{$tm^{xy}_x \sslash sw_{xz} = sw_{xz} \sslash sw_{yz} \sslash tm_x^{xy}$}
\item{$hm^{yz}_y \sslash sw_{xy} = sw_{xy} \sslash sw_{xz} \sslash hm^{yz}_y$}
\item{$sw_{xy} \sslash t\rho^x_u= t\rho^x_u \sslash sw_{uy}$}
\item{$sw_{xy} \sslash h\rho^y_u= h\rho^y_u  \sslash  sw_{xu}$}
\item{$te_x \sslash sw_{xy} = te_x$}
\item{$he_y \sslash sw_{xy}= he_y$}
\end{itemize}
\end{itemize}

Note that in a meta-bicrossed product, {$m^{xy}_z=sw_{xy} \sslash hm^{h_xh_y}_{h_z} \sslash tm^{t_xt_y}_{t_z}$ always defines a meta-monoid with $\Gamma_X=\Gamma(X,X)$}

\section{Some verifications: computer program} \label{sec:Programs}

Using \textit{Mathematica}, it is possible to write a very concise implementation of $\beta$-calculus, and use to carry out the algebraic manipulations that prove Theorem 1 and verify Conjecture 1 on a convincing number of knots and links. We do that in several parts below, with all code included.

\subsection{The Program} We start by loading the \textit{Mathematica} package {\tt KnotTheory`}. This is not strictly necessary, and it is only used for comparison with standard evaluations of the Alexander polynomial:

\vskip 3mm
\noindent\includegraphics[scale=0.18]{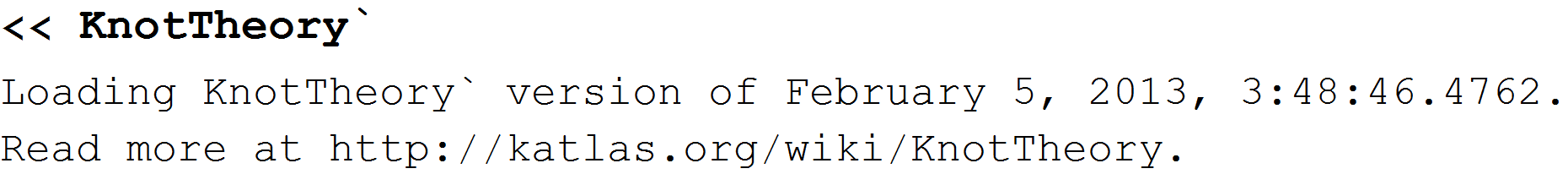}

\vskip 3mm
We then move on to our main program.

\def\annot{{The first part of the program is mostly cosmetic. Its main part is the routine $\beta${\tt Form} used for pretty-printing $\beta$-calculus outputs.}}
\if\jktr y
  \annot
  \par\noindent\imagetop{\includegraphics[scale=0.18]{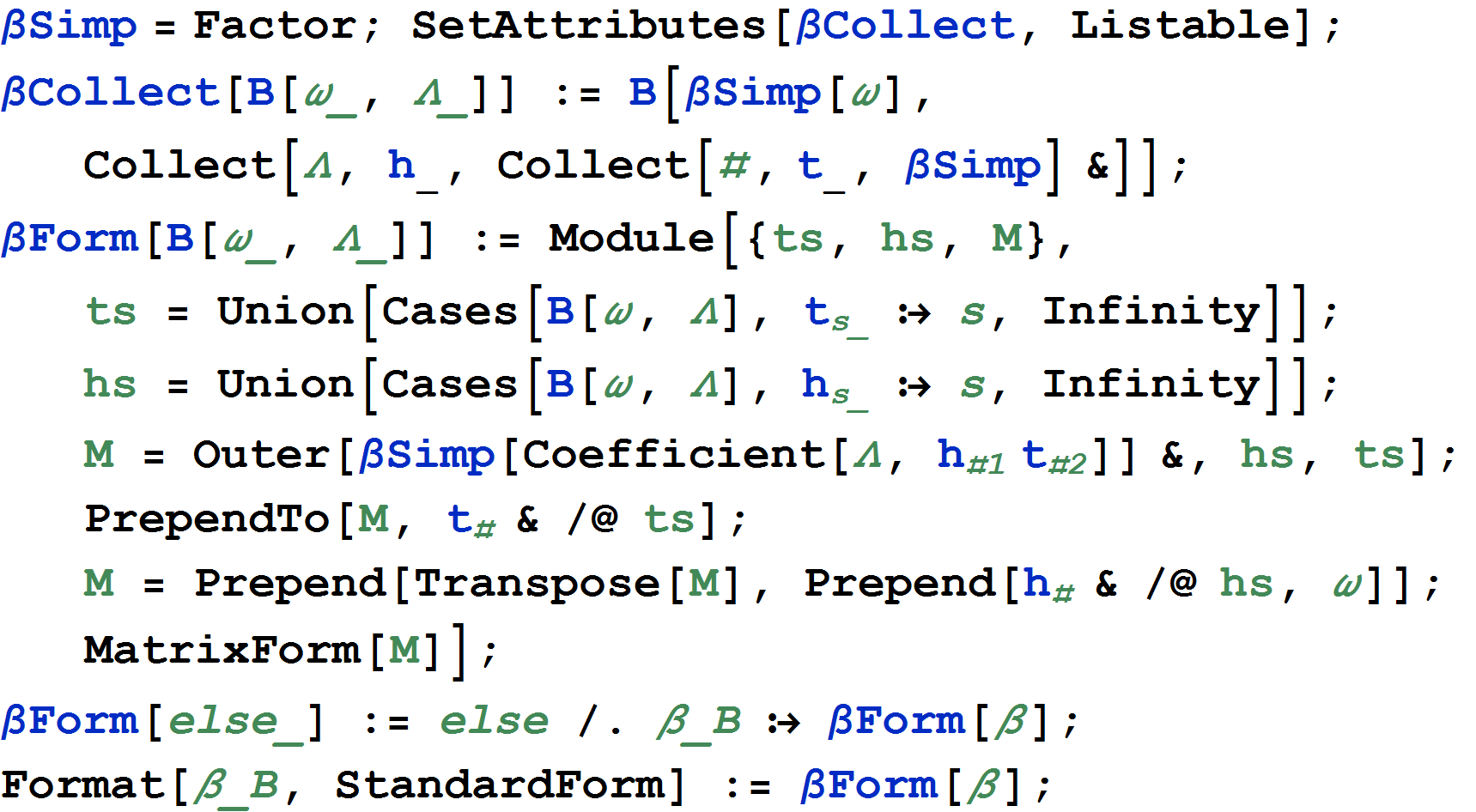}}
\else
  \noindent\imagetop{\includegraphics[scale=0.18]{Initialization.png}}%
  \hfill\raisebox{-3.5mm}{\parbox[t]{1.75in}{\annot}}
\fi

\def\annot{{In the main part of the  program, a $\beta$ matrix is represented as a polynomial in two variables: $\mu = \sum \alpha_{ij} t_i h_j$. This makes some calculations very simple! Selecting the content of column $i$ is achieved by taking a derivative with respect to $h_i$; setting all the $t$'s equal to 1 computes its column sum. The disjoint union of two matrices is simply the sum of their polynomials.}}
\vskip 3mm
\if\jktr y
  \annot
  \par\noindent\imagetop{\includegraphics[scale=0.18]{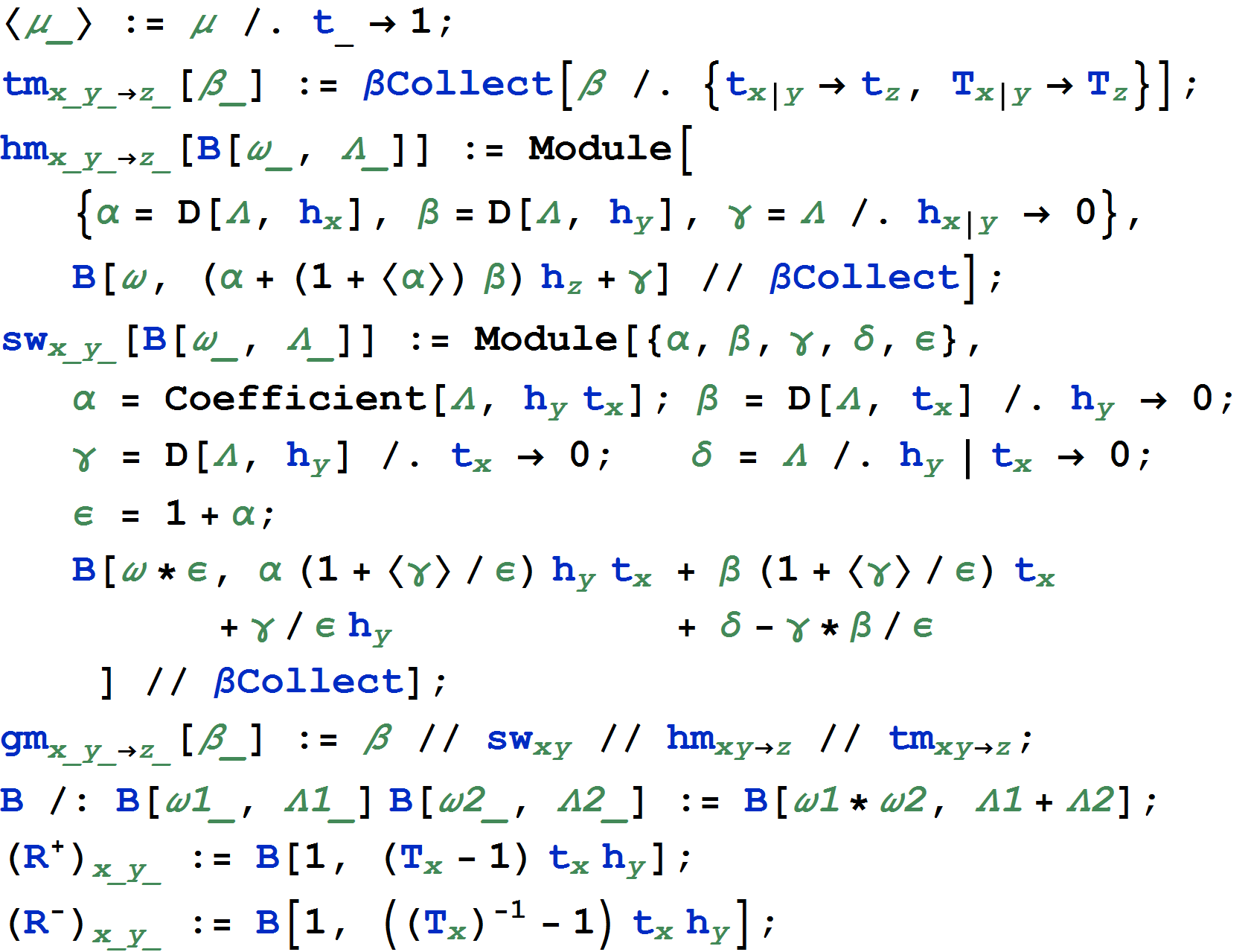}}
\else
  \noindent\imagetop{\includegraphics[scale=0.18]{Program.png}}%
  \hfill\raisebox{-3.5mm}{\parbox[t]{2in}{\annot}}
\fi

\subsection{Proof of Theorem 1}\label{ssec:th1}

To establish Theorem 1 we just need to check that the operations of $\beta$-calculus satisfy the axioms of a meta-bicrossed product listed in Section \ref{ssec:def}. We only bother with the non-obvious axioms, the associativity of $tm$ and of $hm$, and the two swap axioms of Figure~\ref{fig:swap}. Even this we do the lazy way --- we have a computer implementation of the $\beta$-calculus operations. Why not use it to check the relations?

\def\annot{{As a first check, we check the meta-associativity of $tm$ --- we input a generic $4$-tail and $2$-head $\beta$ matrix, let $O_1$ and $O_2$ be the outputs of evaluating $tm^{12}_1\sslash tm^{13}_1$ and $tm^{23}_2\sslash tm^{12}_1$ on $\beta$, and finally we print the logical value of $O_1=O_2$. Nicely, it comes out to be {\tt True}.}}
\if\jktr y
  \annot
  \par\noindent\imagetop{\includegraphics[scale=0.18]{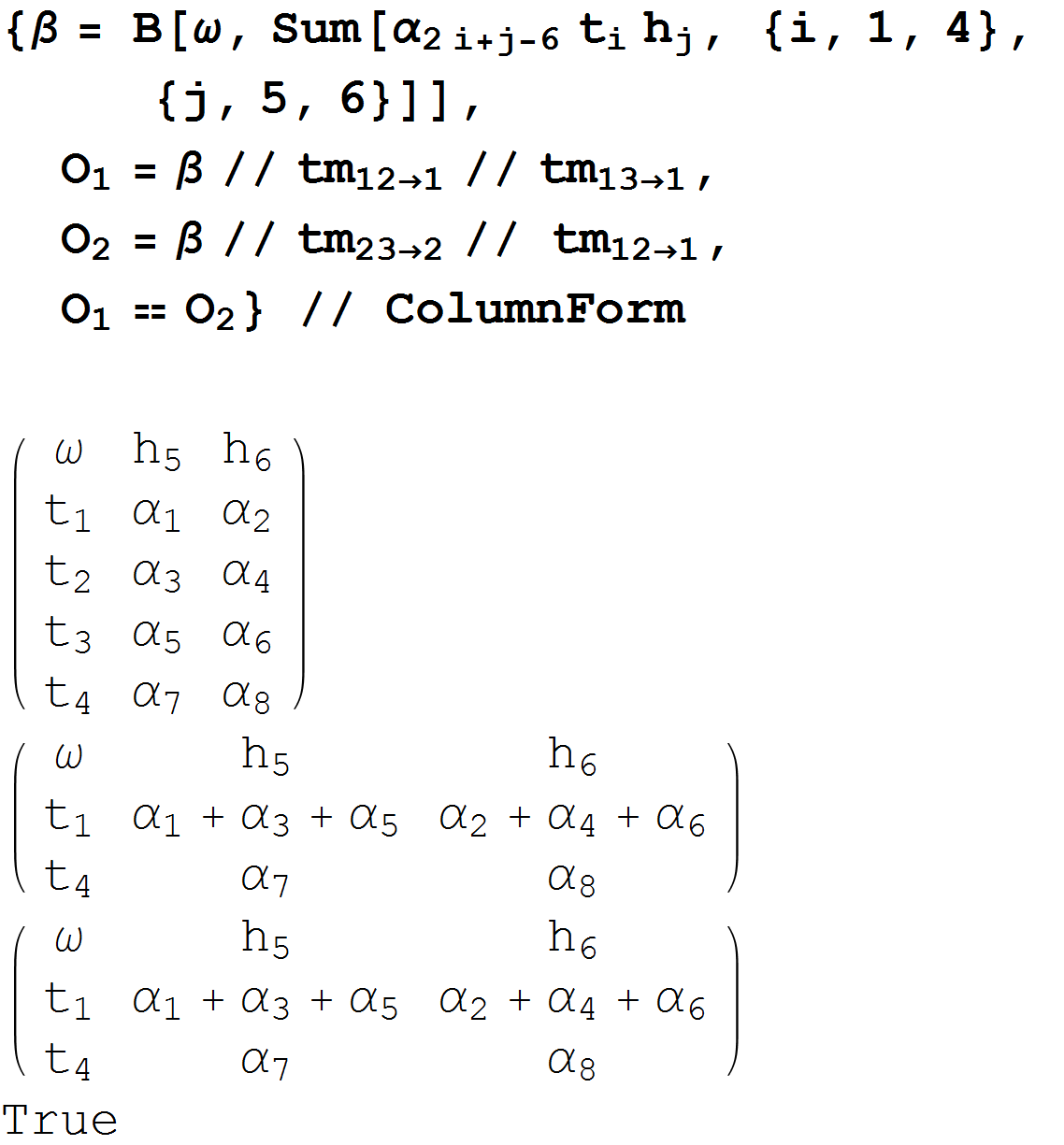}}
\else
  \noindent\imagetop{\includegraphics[scale=0.18]{tm.png}}%
  \hfill\raisebox{-3.5mm}{\parbox[t]{2.5in}{\annot}}
\fi
\vskip 3mm

\def\annot{{We then do the same for $hm$, except we now use a $\beta$ matrix with $2$ tails and $4$ heads, and we suppress the printing of $O_2$. Nicely, the logical value of $O_1=O_2$ is again {\tt True}. (So we didn't lose much by not printing $O_2$). Note that to keep our output from overflowing the width of the page, we have to denote $\alpha_i$ by $\hat{i}$.}}
\if\jktr y
  \annot
  \par\noindent\imagetop{\includegraphics[width=\textwidth]{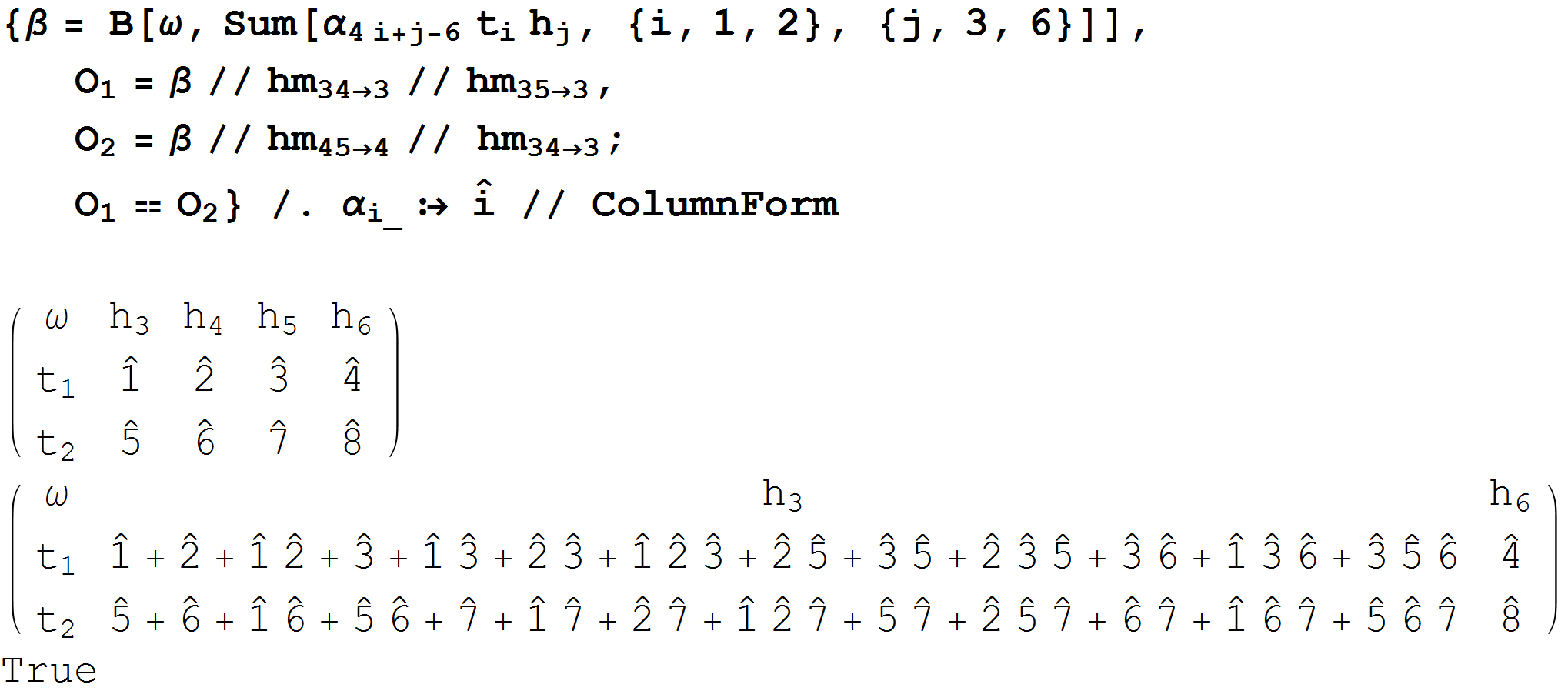}}
\else
  \noindent\imagetop{\includegraphics[scale=0.18]{hm.png}}%
  \hspace{-1.1in}\raisebox{-3.5mm}{\parbox[t]{2.75in}{\annot}}
\fi
\vskip 3mm

\def\annot{{Next come the two swap axioms.}}
\if\jktr y
  \annot
  \par\imagetop{\includegraphics[width=\textwidth]{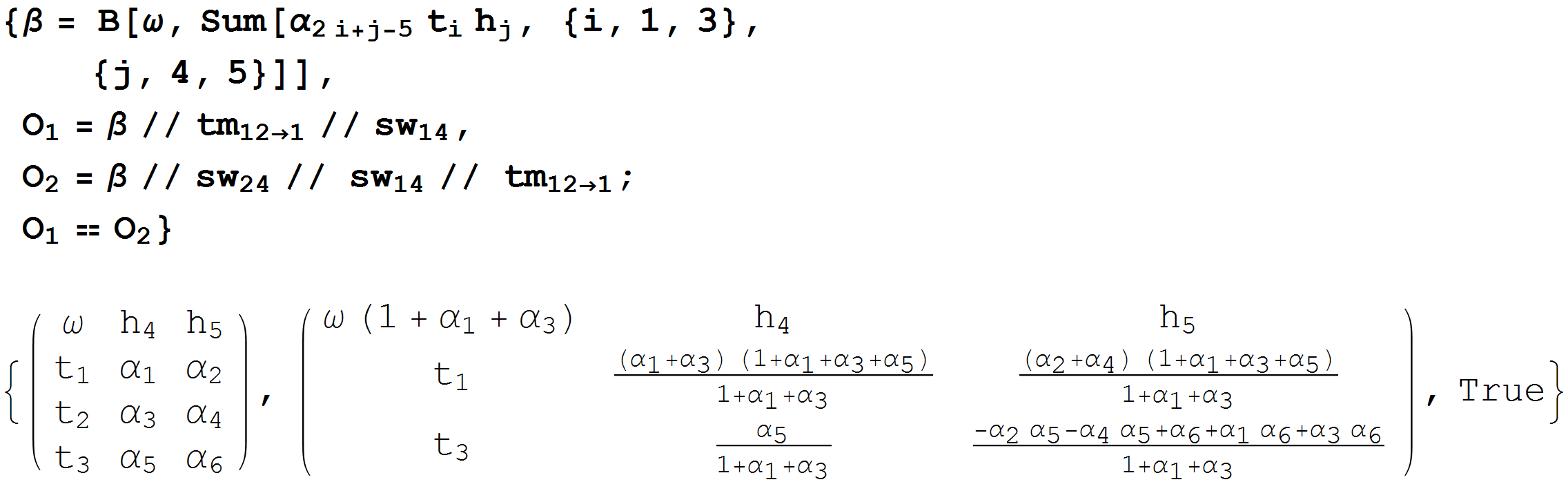}}
\else
  \noindent\imagetop{\includegraphics[scale=0.18]{htt.png}}%
  \hspace{-1.5in}\raisebox{-3.5mm}{\parbox[t]{2.5in}{\annot}}
\fi
\vskip 3mm

\def\annot{{Note that for the second swap axiom, some algebraic simplification must take place, using the routine {\tt $\beta$Collect}.}}
\if\jktr y
  \annot
  \par\noindent\imagetop{\includegraphics[width=\textwidth]{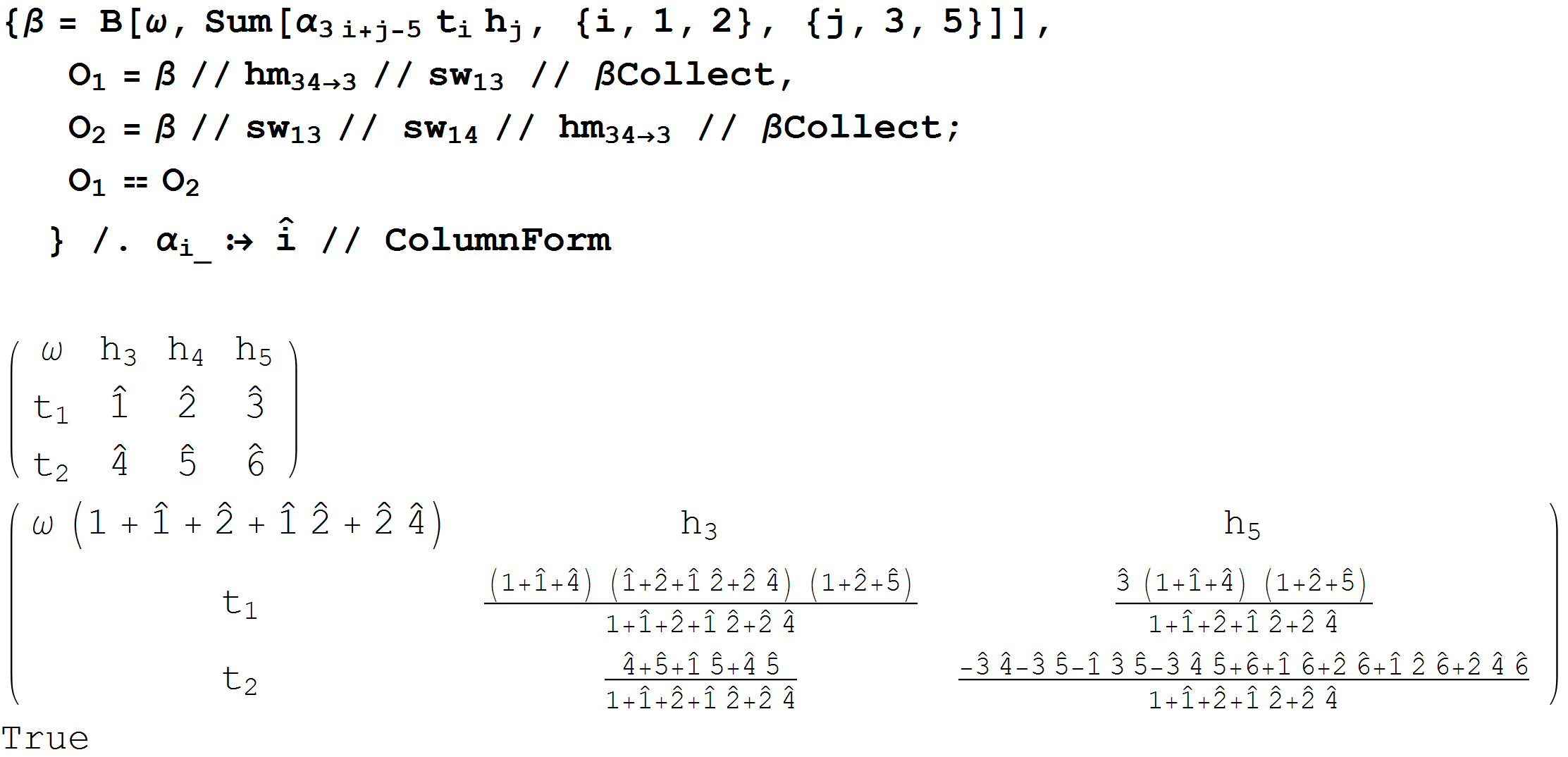}}
\else
  \noindent\imagetop{\includegraphics[scale=0.18]{hht.png}}%
  \hspace{-1.5in}\raisebox{-3.5mm}{\parbox[t]{2.5in}{\annot}}
\fi
\vskip 3mm

\def\annot{{Just for completeness, we verify the third Reidemeister move once again.}}
\if\jktr y
  \annot
  \par\noindent\imagetop{\includegraphics[scale=0.18]{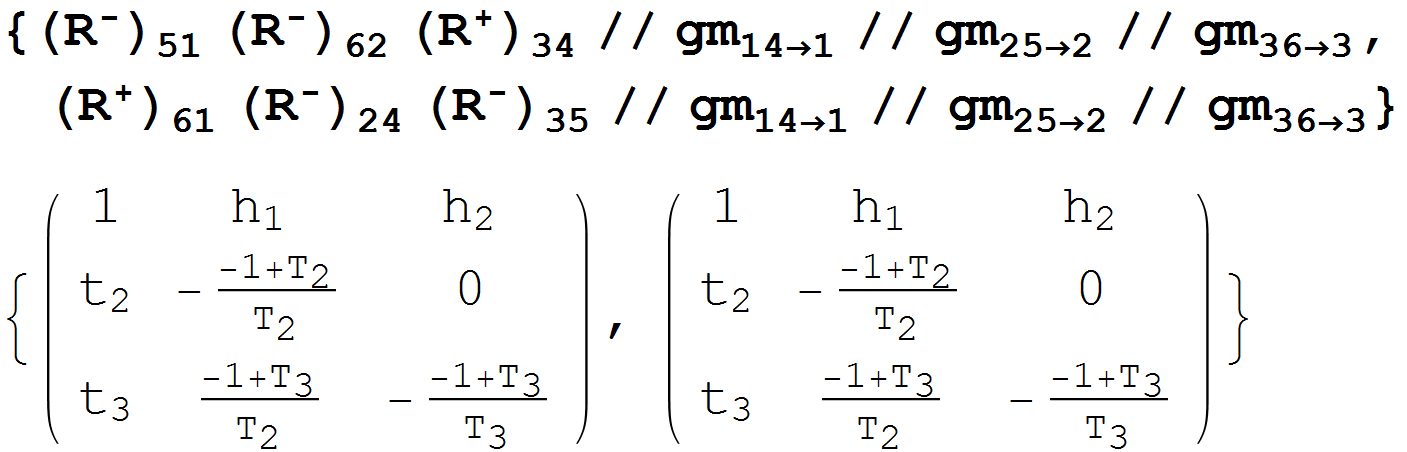}}
\else
  \noindent\imagetop{\includegraphics[scale=0.18]{R3.png}}%
  \hfill\raisebox{-3.5mm}{\parbox[t]{2.5in}{\annot}}
\fi
\vskip 3mm

\subsection{Testing Conjecture~\ref{conj:Alexander}} \label{ssec:KnotsAndLinks}
Our next task is to carry out some computations for knots and links in support of Conjecture~\ref{conj:Alexander}. As our first demonstration, we compute $Z^\beta(8_{17})$ in several steps. The first step is to generate the invariant of the tangle consisting of the disjoint union of $8$ crossings, labeled as the crossings of $8_{17}$ are labeled but not yet connected to each other:

\noindent\imagetop{\includegraphics[scale=0.18]{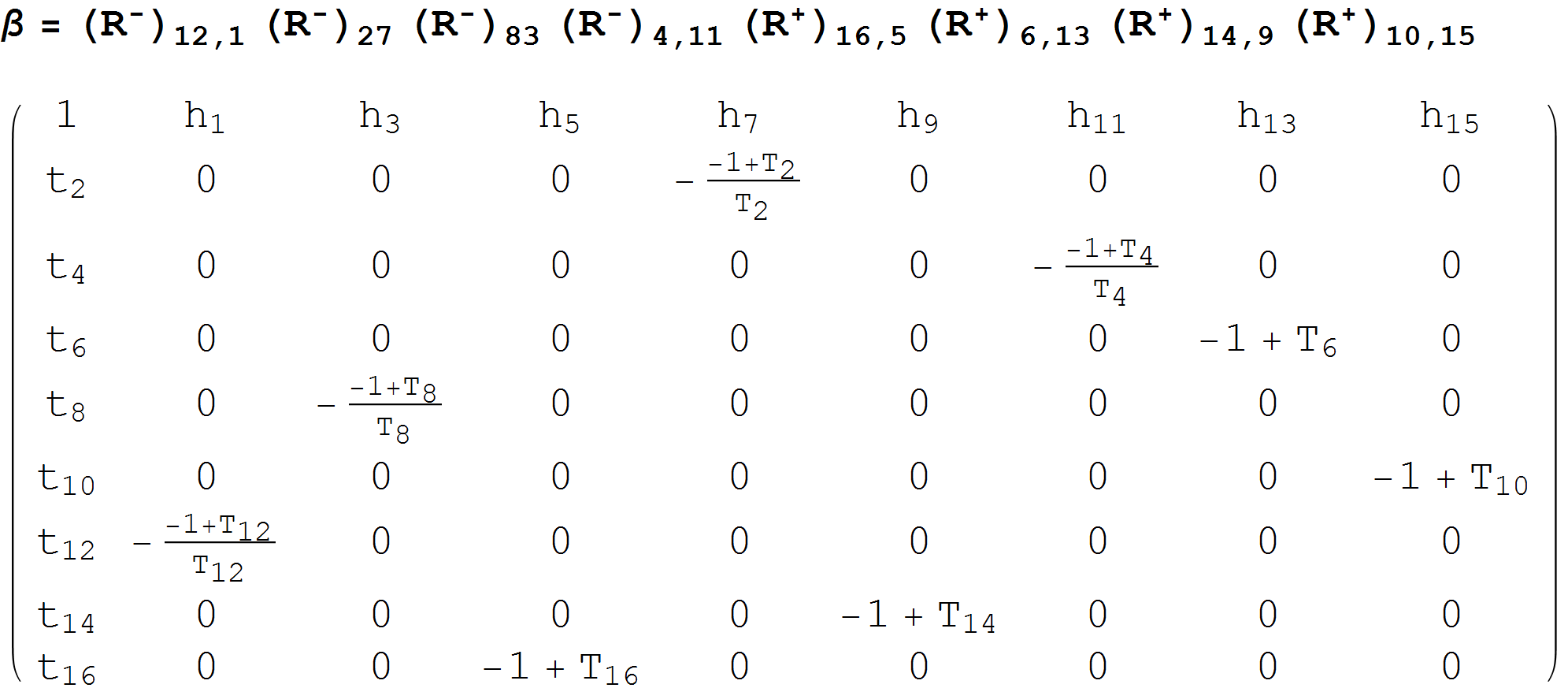}}%
\vskip 3mm

Next, we {\em partially} concatenate the strands of these 8 crossings to each other, making only 9 of the required 15 connections. The result is 3-component tangle that approximates $8_{17}$, and a chance to see what an intermediate step of the computation looks like:

\if\jktr y
  \noindent\imagetop{\includegraphics[width=\textwidth]{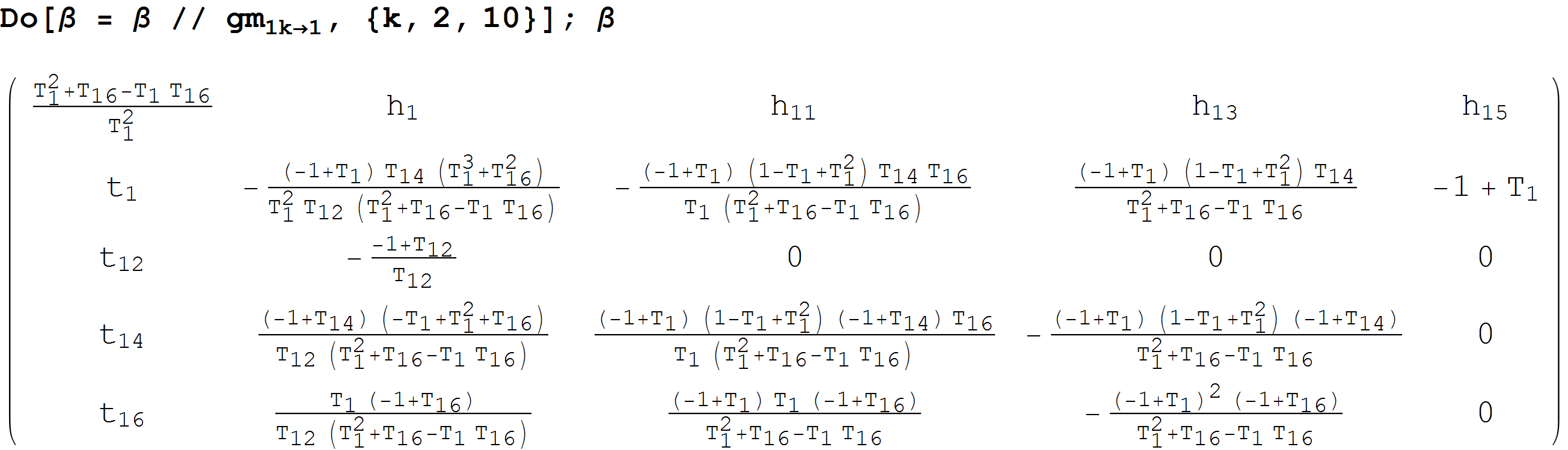}}%
\else
  \noindent\imagetop{\includegraphics[scale=0.18]{8_17-2.png}}%
\fi
\vskip 3mm

\def\annot{{We then complete the sewing together of $8_{17}$, obtaining $Z^\beta(8_{17})$. Note that the ``matrix part'' of the invariant is completely suppressed by our printing routine, because it is $0$.}}
\if\jktr y
  \annot
  \par\noindent\imagetop{\includegraphics[scale=0.18]{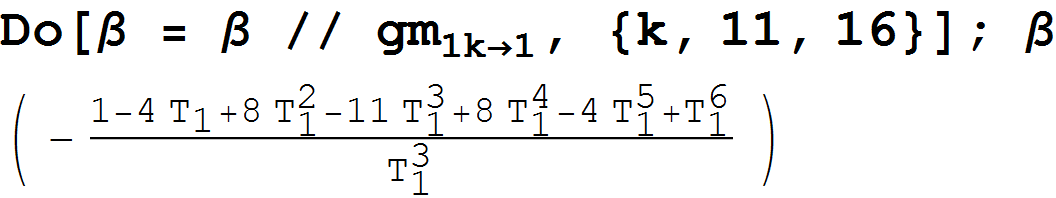}}
\else
  \noindent\imagetop{\includegraphics[scale=0.18]{8_17-3.png}}%
  \hfill\raisebox{-3.5mm}{\parbox[t]{3.5in}{\annot}}
\fi
\vskip 3mm

\def\annot{{For completeness, we compare with the pre-computed value of the Alexander polynomial, as known to {\tt KnotTheory`}. As can be fairly expected, it differs from the computed value of $Z^\beta(8_{17})$ by a unit.}}
\if\jktr y
  \annot
  \par\noindent\imagetop{\includegraphics[scale=0.18]{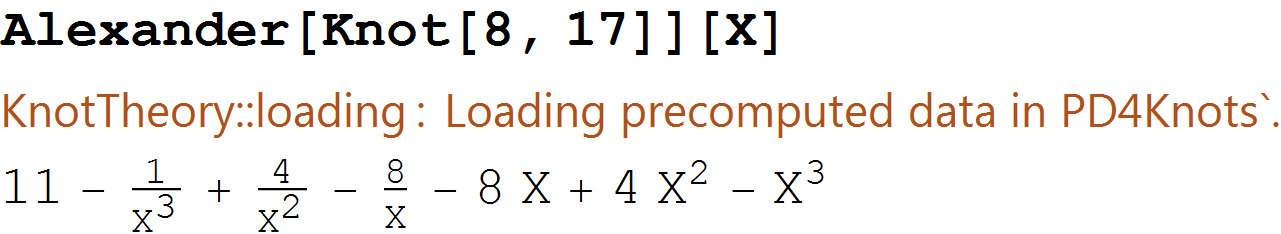}}
\else
  \noindent\imagetop{\includegraphics[scale=0.18]{8_17-4.png}}%
  \hfill\raisebox{-3.5mm}{\parbox[t]{3.1in}{\annot}}
\fi
\vskip 3mm

\def\annot{{We next make it systematic by writing a short program that compute $Z^\beta$ of an arbitrary input link.}}
\if\jktr y
  \annot
  \par\noindent\imagetop{\includegraphics[scale=0.18]{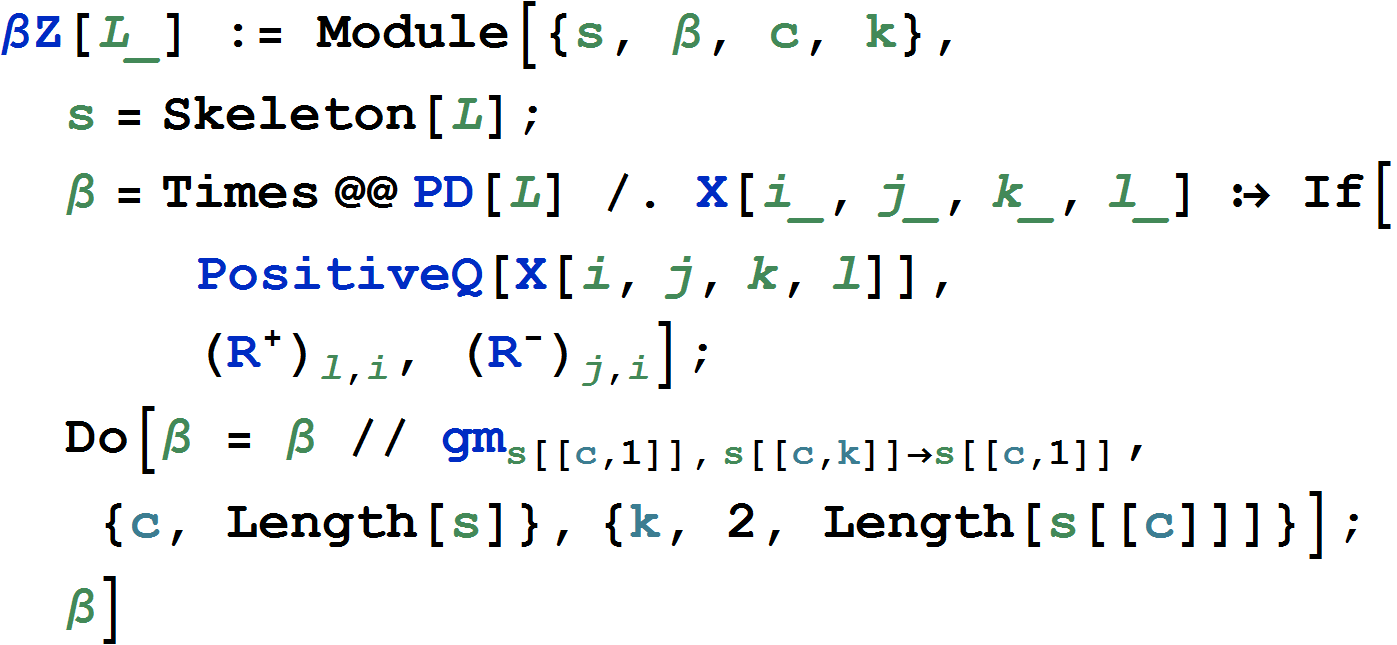}}
\else
  \noindent\imagetop{\includegraphics[scale=0.18]{betaZ.png}}%
  \hfill\raisebox{-3.5mm}{\parbox[t]{2.5in}{\annot}}
\fi
\vskip 3mm

\def\annot{{We verify that for all knots with up to 8 crossings, the ratio of $Z^\beta$ and the Alexander polynomial is always a unit. At home we've verified the same thing for all knots with up to 11 crossings.}}
\if\jktr y
  \annot
  \par\noindent\imagetop{\includegraphics[scale=0.18]{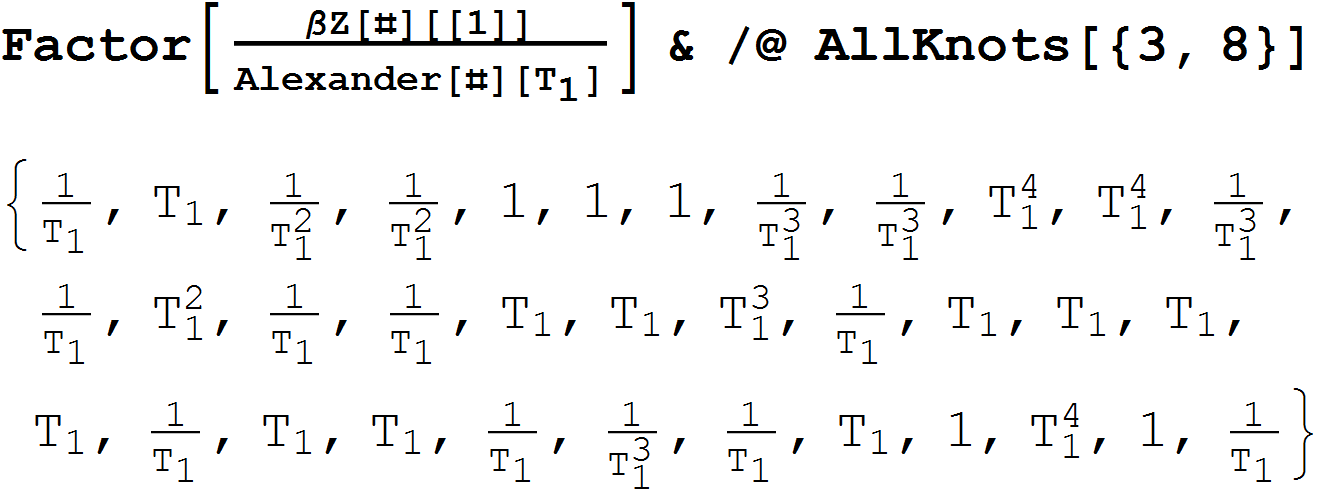}}
\else
  \noindent\imagetop{\includegraphics[scale=0.18]{TestAllKnots.png}}%
  \hfill\raisebox{-3.5mm}{\parbox[t]{2.5in}{\annot}}
\fi
\vskip 3mm

\def\annot{{Next is the program for extracting the multi-variable Alexander polynomial from the information in $Z^\beta$.}}
\if\jktr y
  \annot
  \par\noindent\imagetop{\includegraphics[scale=0.18]{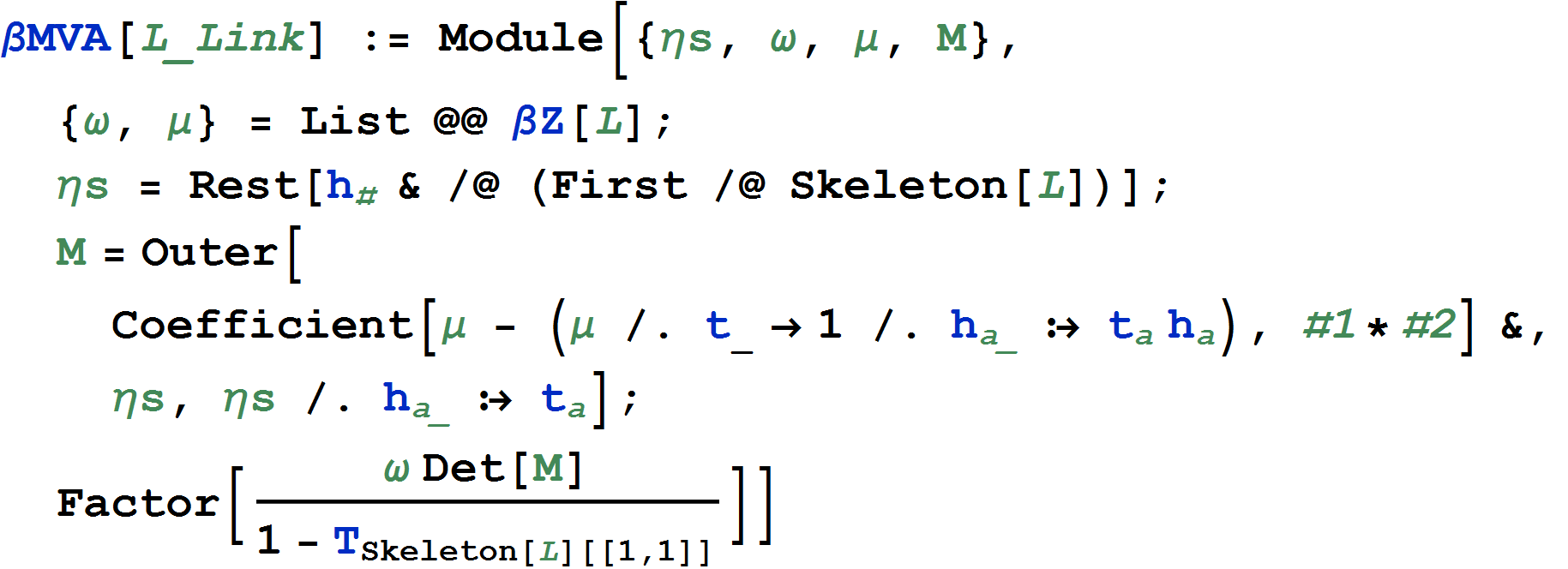}}
\else
  \noindent\imagetop{\includegraphics[scale=0.18]{betaMVA.png}}%
  \hfill\raisebox{-3.5mm}{\parbox[t]{1.7in}{\annot}}
\fi
\vskip 3mm

\def\annot{{It works for the Borromean rings!}}
\if\jktr y
  \annot
  \par\noindent\imagetop{\includegraphics[scale=0.18]{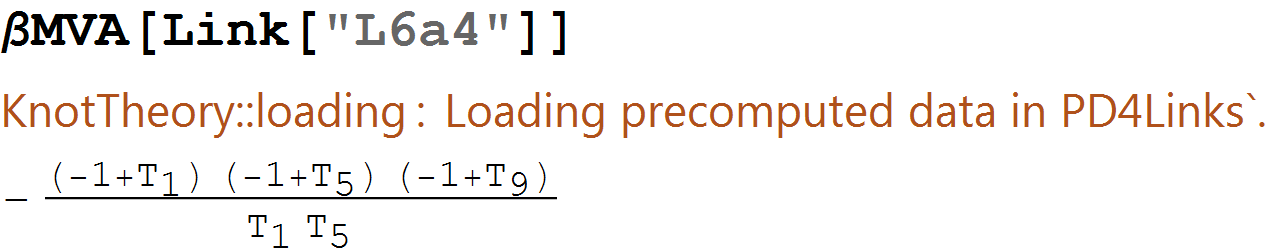}}
\else
  \noindent\imagetop{\includegraphics[scale=0.18]{BorromeanMVA.png}}%
  \hfill\raisebox{-3.5mm}{\parbox[t]{1.5in}{\annot}}
\fi
\vskip 3mm

\def\annot{{And also for all links with up to 7 crossings. At home we've verified the same for all links with up to 11 crossings.}}
\if\jktr y
  \annot
  \par\noindent\imagetop{\includegraphics[scale=0.18]{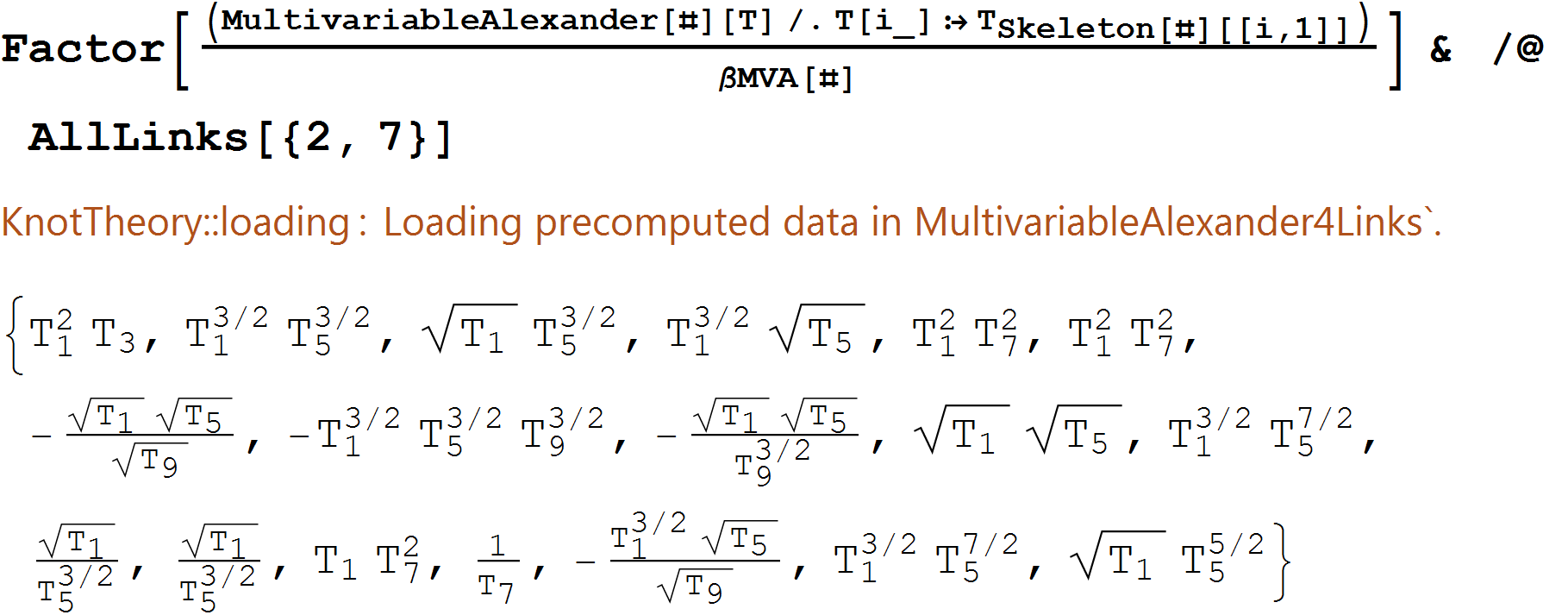}}
\else
  \noindent\imagetop{\includegraphics[scale=0.18]{TestAllLinks.png}}%
  \hfill\raisebox{-3.5mm}{\parbox[t]{1.5in}{\annot}}
\fi
\vskip 3mm

\section{Acknowledgement}
\if\jktr y
  This work was partially supported by NSERC grant RGPIN 262178 and partially pursued at the Newton Institute in Cambridge, UK.
\fi
We wish to thank Iva Halacheva for comments and suggestions.

\end{document}